\documentclass[10pt,a4paper,oneside]{article}

\usepackage{amsmath}
\usepackage{amsfonts}
\usepackage{amssymb}
\usepackage{amsthm}
\usepackage{latexsym}
\usepackage[active]{srcltx}

\usepackage{enumerate}
\usepackage{color}

 \sloppypar

\newtheorem{theorem}{Theorem}[section]

\newtheorem{proposition}[theorem]{Proposition}

\renewcommand{\H}{{\mathcal H}}
\def\al{\alpha}
\def\be{\beta}

\def\de{\delta}

\def\et{\eta}

\def\ga{\gamma}
\def\GA{\Gamma}

\def\la{\lambda}

\def\om{\omega}
\def\va{\varphi}

\def\g{\mathfrak{g}}
\def\a{\mathfrak{a}}

\def\q{\mathfrak{q}}
\def\p{\mathfrak{p}}
\def\n{\mathfrak{n}}

\def\s{\mathfrak{s}}
\def\t{\mathfrak{t}}
\def\z{\mathfrak{z}}
\def\u{\mathfrak{u}}
\def\r{\mathfrak{r}}

\def\l{\mathfrak l}

\def\ga{\gamma}
\def\la{\lambda}

\def\si{\sigma}
\def\om{\omega}

\def\ch{\chi}

\def\ps{\psi}

\def\ol#1{\overline{#1}}
\def\nn{\nonumber}

\def\R{{\mathbb R}}
\def\C{{\mathbb C}}
\def\N{{\mathbb N}}

\def\Id{{\mathbb I}}

\def\B{{\mathcal B}}

\def\F{{\mathcal F}}

\def\H{{\mathcal H}}

\def\K{{\mathcal K}}

\def\M{{\mathcal M}}

\def\S{{\mathcal S}}

\def\ad{{\text ad}}
\def\Ad{{\text Ad}}

\def\iy{\infty}

\def\wh{\widehat}

\def\ol#1{\overline{#1}}

\def\hb#1{\hbox{#1}}
\def\val#1{\vert #1\vert}

\def\no#1{\Vert #1\Vert }

\def\span#1{\text{span}\{#1\}}
\def\exp{\text{exp}}

\def\res#1{_{\vert #1}}
\def\inv{^{-1}}

\def\es{\emptyset}

\def\hb #1{\hbox{#1}}

\def\hb#1{\hbox{#1}}
\def\val#1{\vert #1\vert}

\def\noop#1{\Vert #1\Vert_{\rm op}}

\def\dim#1{\hb{dim}(#1)}

\def\ind{{\rm ind} }

\def\L1#1{L^1(#1)}

\def\L#1#2{L^{#1}(#2)}
\def\l#1#2{L^{#1}(#2)}

\def\ci{\circ}

\def\lef({\left(}
\def\rig){\right)}

\def\ot{\otimes}

 \sloppypar

\usepackage{a4wide}

\newtheorem{definition}[theorem]{Definition}

\newtheorem{conditions}[theorem]{Conditions}
\newtheorem{lemma}[theorem]{Lemma}

\newcommand{\klgl}{\leq}

\newcommand{\quer}{\overline}

\def\ad{\mathrm{\, ad\, }}
\def\Ad{\mathrm{\, Ad \,}}

\def\a{\mathfrak a}

\def\u{\mathfrak u}
\def\z{\mathfrak z}
\def\n{\mathfrak n}

\begin{document}

\title{The $C^*$-algebras of connected real two-step nilpotent Lie groups}
\author{Janne-Kathrin G\"unther and Jean Ludwig}
\maketitle

~\\
\begin{abstract}
Using the operator valued Fourier transform, the 
$C^*$-algebras of connected real two-step nilpotent Lie groups are  
characterized as algebras of operator fields defined over their spectra. In 
particular, it is shown by explicit computations, that the Fourier transform 
of such \mbox{$C^*$-algebras} fulfills the norm controlled dual limit property.
\end{abstract}
~\\

\section{Introduction}
\label{intro}
In this article, the structure of the $C^*$-algebras of two-step nilpotent Lie groups will be analyzed. In order to be able to understand these $C^*$-algebras, the Fourier transform is an important tool. The Fourier transform $\F(a)=\hat{a}$ of an element $a$ of a $C^*$-algebra $A$ is defined in the following way: One chooses for every $\gamma$ in $\wh{A}$, the spectrum of $A$, a representation $(\pi_{\gamma},\H_{\gamma})$ in the equivalence class of $\gamma$ and defines
$$\F(a)(\gamma):=\pi_{\gamma}(a) \in \H_{\gamma}~~~\forall ~\gamma \in \wh{A}.$$
Then $\F(a)$ is contained in the algebra of all bounded operator fields over $\wh{A}$
$$l^{\infty}(\wh{A})= \Big\{ \phi= \big( \phi(\pi_{\ga}) \in \B(\H_{\ga})\big)_{\ga \in \wh{A}}~|~ \| \phi \|_{\infty}:= \sup \limits_{\ga \in \wh{A}} \|\phi(\pi_{\ga})\|_{op}< \infty \Big\}$$
and the mapping $$\F: A \to l^{\infty}(\wh{A}),~ a \mapsto \hat{a}$$
is an isometric $*$-homomorphism. \\
~\\
The structure of the $C^*$-algebras is already known for certain classes of Lie groups: The $C^*$-algebras of the Heisenberg and the thread-like Lie groups have been charac\-terized in \cite{ludwig-turowska} and the $C^*$-algebras of the $ax+b$-like groups in \cite{lin-ludwig}. Furthermore, the $C^*$-algebras of the 5-dimensional nilpotent Lie groups have been determined in \cite{hedi} and H.Regeiba analyzed the $C^*$-algebras of all 6-dimensional nilpotent Lie groups in his doctoral thesis (see \cite{hedidr}). \\
The methods in this paper will partly be similar, but more complex, to the one used for the characterization of the $C^*$-algebra of the Heisenberg Lie group (see \cite{ludwig-turowska}), which is also two-step nilpotent and thus serves as an example. \\
It will be shown that the $C^*$-algebras of two-step nilpotent Lie groups $G$ are charac\-terized by the following conditions. The same conditions hold true for all 5- and 6-dimensional nilpotent Lie groups (see \cite{hedi}), for the Heisenberg Lie groups and the thread-like Lie groups (see \cite{ludwig-turowska}).

\pagebreak

\begin{conditions}
~\\
\begin{enumerate}
\item Stratification of the spectrum: 
\begin{enumerate}
\item A finite increasing family $ S_0\subset S_1\subset\ldots\subset S_r=\wh{C^*(G)} \cong \wh{G} $ of closed
subsets of the spectrum $ \wh{C^*(G)} \cong \wh{G} $ of $ C^*(G)$ or respectively $G$ will be constructed in such a way that for $ i \in \{1,\cdots, r\} $ the subsets $\GA_0=S_0$ and
$ \GA_i:=S_i\setminus S_{i-1}$ are Hausdorff in their relative topologies and such that $S_0$ consists of all the characters of $C^*(G)$ or $G$, respectively. 
\item For every $ i\in \{0,\cdots, r\} $ a Hilbert space $ \H_i $ and for every $ \ga\in \GA_i $ a concrete realization $ (\pi_\ga,\H_i) $ of $ \ga $ on the Hilbert space $ \H_i $ will be defined.
\end{enumerate}
\item CCR $C^*$-algebra: \\
It will be shown that $C^*(G)$ is a separable CCR (or liminal) $C^*$-algebra, i.e. a separable \mbox{$C^*$-algebra} such that the image of every irreducible representation
$ (\pi,\H)  $ of $ C^*(G) $ is contained in the algebra of compact operators $ \K(\H) $ (which implies that the image equals $\K(\H)$).
\item Changing of layers: \\
Let $a \in C^*(G)$. 
\begin{enumerate}
\item  It will be proved that the mappings $ \ga \mapsto \F (a)(\ga) $ are norm continuous on the different sets $ \GA_i $.
\item   For any  $ i \in \{0,\cdots, r\} $  and for any
converging sequence contained in $ \GA_i $ with limit set outside $ \GA_i $ (thus in $S_{i-1}$), there will be constructed a properly 
converging subsequence $\ol\ga=(\ga_k)_{k\in\N} $ (i.e. the subsequences of $\ol\ga$ have all the same limit set - see Definition \ref{propconv}), as well as a constant $ C>0 $ and for every $ k\in\N $ an involutive linear mapping $ \tilde{\nu}_k=\tilde{\nu}_{\ol\ga,k}: CB(S_{i-1})\to \B(\H_i)$, {which is  bounded by $ C\|\cdot\|_{S_{i-1}} $}, such that
\begin{eqnarray*}
&\lim \limits_{k\to\iy }\big\|
\F (a)(\ga_k)-\tilde{\nu}_{k}\big(\F (a)\res{S_{i-1}}\big)\big\|_{op}=0. 
\end{eqnarray*}
Here $CB(S_{i-1})$ is the $ * $-algebra of all the uniformly
bounded fields of operators \linebreak
$ \big(\ps(\ga)\in \B(\H_j) \big)_{\ga \in \GA_j, j=0,\cdots, i-1}$, which are
operator norm continuous on the subsets $ \GA_j$ for every $  j\in\{0,\cdots, i-1\} $, provided with the infinity-norm 
\begin{eqnarray*}
\no{\va}_{S_{i-1}}:=\sup_{\ga\in S_{i-1}}\noop{\va(\ga)}. 
\end{eqnarray*}
 \end{enumerate}
 \end{enumerate}
\end{conditions}

These properties characterize the structure of  $C^*(G)$ (see \cite{hedi}, Theorem 3.5). A $C^*$-algebra fulfilling these conditions is called a $C^*$-algebra with "norm controlled dual limits". \\ 
The main work of this article consists in the proof of Property 3(b) and in particu\-lar in the construction of the mappings $(\tilde{\nu}_k)_k$.

\section{Preliminaries} \label{pre}
\subsection{Two-step nilpotent Lie groups}
Let $\g$ be a real Lie algebra which is nilpotent of step two. This means that 
\begin{eqnarray*}
 [\g,\g]:=span \big\{[X,Y]|~ X,Y \in \g \big\}
 \end{eqnarray*}
is contained in the center of $\g $. \\
Fix a scalar product $ \langle \cdot,\cdot\rangle $ on $ \g $ and take on $ \g $ the Campbell-Baker-Hausdorff multiplication
\begin{eqnarray*}
 u\cdot v &= &u+v+\frac{1}{2}[u,v]~~~\forall~ u,v\in\g.
\end{eqnarray*}

This gives the simply connected connected Lie group $G=(\g,\cdot)$ with Lie algebra $\g $. The exponential mapping $\exp: \g\to G=(\g,\cdot)$ is in this case the identity mapping. \\

The Haar measure of this group is a Lebesgue measure which is denoted by $dx$. \\
Then, the $C^*$-algebra of $G$ is defined as the completion of the convolution algebra \mbox{$\l1{G,dx}=\l1G $} with respect to the $C^*$-norm of $\l1{G,dx}$, i.e. 
$$C^*(G)~:= ~\quer{L^1(G,dx)}^{\| \cdot \|_{C^*(G)}}~~~\text{with}~~~\|f \|_{C^*(G)}~:=~\sup \limits_{\pi \in \wh{G}}\| \pi(f) \|_{op}$$
and a well-known result, that can be found in \cite{dix}, states that the spectrum of $C^*(G)$ coincides with the spectrum of $G$: $$\wh{C^*(G)}=\wh{G}.$$
~\\
Now, for a linear functional $\ell$ of $\g$, consider the skew-bilinear form 
\begin{eqnarray*}
 B_{\ell}(X,Y):=\langle{\ell},{[X,Y]}\rangle
 \end{eqnarray*}
on $\g $. Moreover, let 
\begin{eqnarray*}
 \g(\ell):=\big\{X\in\g|~\langle{\ell},{[X,\g]\rangle=\{0\}} \big\}
 \end{eqnarray*}
be the radical of $B_\ell $ and the stabilizer of the linear functional $\ell$. Then, as $\g$ is two-step nilpotent, $[\g,\g]\subset \g(\ell)$ and thus $\g(\ell)$ is an ideal of $\g$. \\
\begin{definition}
~\\
\rm
A subalgebra $\p $ of $\g $, that is \textit{subordinated} to $\ell $ (i.e. that fulfills $\langle{\ell},{[\p,\p]}\rangle=\{0\} $) and that has the dimension
\begin{eqnarray*}
 \dim \p=\frac{1}{2} \big(\dim \g+\dim{\g(\ell)} \big),
 \end{eqnarray*}
 which means that $\p $ is maximal isotropic for $B_\ell $, is called a \textit{polarization} in $\ell$. 
\end{definition}
Again since $\g $ is nilpotent of step two, every maximal isotropic subspace $\p $ of $\g $ for $B_\ell $ containing $[\g,\g] $ is a polarization at $\ell $. \\
Now, if $\p\subset \g $ is any subalgebra of $\g $ which is subordinated to $\ell $, the linear functional $\ell $ defines a unitary character $\chi_{\ell}$ of $P:=\exp(\p)$:
$$\ch_{\ell}(x)~:=~e^{-2  \pi i \langle{\ell},{\log(x)}\rangle}=~e^{-2  \pi i \langle{\ell},{x}\rangle}~~~\forall~x\in P.$$

\subsection{Induced representations}
The induced representation $\si_{\ell,\p}=\ind_ P ^G\ch_\ell $ for a polarization $\p $ in $\ell $ and $P:=exp(\p)$ can be described in the following way: 

Since $\p $ contains $[\g,\g] $ and even the center $\z $ of $\g $, one can write $\g=\s\oplus\p$ and $\p=\t\oplus\z $ for two subspaces $\t$ and $\s$ of $\g $. 
The quotient space $G/P $ is then homeomorphic to $ \s$ and the Lebesgue measure $d s $ on $ \s$ defines an invariant Borel measure $d\dot g $ on $G/P $. The group $G $ acts by the left translation $\si_{\ell,\p} $ on the Hilbert space
\begin{eqnarray*}
 L^2(G/P,\chi_{\ell})&:=& \bigg\{\xi:G\to \C|~ \xi \text{ measurable},~\xi(gp)=\ol{\ch_\ell(p)}\xi(g)~\forall~ g\in G~\forall~ p\in P,\\
&&~~~\no{\xi}_2^2:=\int \limits_{G/P}\val{\xi(g)}^2~d\dot g<\iy \bigg\}.
 \end{eqnarray*}

Now, if one uses the coordinates $G=\s\cdot \p$, one can identify the Hilbert spaces $L^2(G/P,\chi_{\ell})$ and $L^2(\s)=L^2(\s,ds) $: \\
Let $U_{\ell}:L^2(\s,ds) \to  L^2(G/P,\chi_{\ell})$ be defined by
\begin{eqnarray}\label{idltwtw}
 \nn U_\ell(\va)(S\cdot Y) ~:= ~\ch_{\ell}(-Y)\va(S)~~~\forall ~Y\in\p~~\forall ~S\in\s~~\forall ~\va\in L^2(\s). 
\end{eqnarray}
Then, $U_{\ell}$ is a unitary operator and one can transform the representation $\si_{\ell,\p} $ into a representation $\pi_{\ell,\p} $ on the space $ L^2(\s)$:
\begin{eqnarray}\label{indident}
 \pi_{\ell,\p}~:=~U_\ell^*\ci \si_{\ell,\p}\ci U_\ell.
 \end{eqnarray}
Furthermore, one can express the representation $\si_{\ell,\p}$ in the following way:
\begin{eqnarray*}
 \nn \si_{\ell,\p}(S\cdot Y)\xi(R) &= &\xi \big(Y\inv S\inv R \big)\\
 \nn  &=& \xi \Big((R-S) \cdot \Big(-Y+\frac{1}{2}[R,S]-\frac{1}{2}[R-S,Y] \Big)\Big)\\
 \nn  &=&  e^{2 \pi i\langle{\ell},{-Y+\frac{1}{2}[R,S]-\frac{1}{2}[R-S,Y]}\rangle}\xi(R-S)\\
 \nn&&~\forall ~R,S \in\s~~\forall ~Y\in\p~~\forall ~\xi\in L^2(G/P,\chi_{\ell}).
\end{eqnarray*}
Hence
\begin{eqnarray}\label{exppi}
  \pi_{\ell,\p}(S \cdot Y)\va(R)&=&e^{2\pi i\langle{\ell},{-Y+\frac{1}{2}[R,S]-\frac{1}{2}[R-S,Y]}\rangle}\va(R-S)\\
  \nn&&~\forall~ R,S \in\s~~\forall~ Y\in\p~~\forall~ \va\in L^2(\s).  
\end{eqnarray} 

\subsection{Orbit method}
By the \textbf{Kirillov theory} (see \cite{cor-green}, Section 2.2), for every representation class $\gamma \in \wh{G}$, there exists an element $\ell \in \g^*$ and a polarization $\p$ of $\ell$ in $\g$ such that $\gamma=\big[ind_P^G \chi_{\ell}\big]$, whereat $P:=exp (\p)$. Moreover, if $\ell, \ell' \in \g^*$ are located in the same coadjoint orbit $O \in \g^*/G$ and $\p$ and $\p'$ are polariza\-tions in $\ell$ and $\ell'$, respectively, the induced representations $ind_{P}^G \chi_{\ell}$ and $ind_{P'}^G \chi_{\ell'}$ are equivalent and thus, the Kirillov map which goes from the coadjoint orbit space $\g^*/G $ to the spectrum $\wh G $ of $G $
$$K: \g^*/G \to \wh{G},~ Ad^*(G)\ell \mapsto \big[ind_P^G \chi_{\ell} \big]$$
is a homeomorphism (see \cite{brown} or \cite{leptin-ludwig}, Chapter 3). Therefore,
$$\g^*/G \cong \wh{G}$$
as topological spaces. \\
For every $\ell\in\g^*$ and $x\in G=(\g,\cdot)$  
\begin{eqnarray*}
\Ad^*(x)\ell=\big(\Id_{\g^*}+\ad^*(x) \big)\ell \in \ell+\g(\ell)^\perp.
 \end{eqnarray*}
Hence, as $\ad^*(\g)\ell= \g(\ell)^\perp$,
\begin{eqnarray}\label{desor}
 O_\ell:=\Ad^*(G)\ell=\ell+\g(\ell)^\perp.
 \end{eqnarray}
 
\begin{definition}\label{propconv}
~\\
\rm  Let $T $ be a second countable topological space and suppose that $T $ is not Hausdorff, which means that converging sequences can have many limit points. Denote by $L ((t_k)_k ) $ the collection of all the limit points of a sequence $(t_k)_k $ in $T $. A sequence $(t_k)_k $ is called \textit{properly converging}, if $(t_k)_k$ has limit points and if every subsequence of $(t_k)_k$ has the same limit set as $(t_k)_k$. \\
It is well known that every converging sequence in $T $ admits a properly converging subsequence.  
\end{definition}

Now, let $(\pi_k)_k\subset \wh G $ be a properly converging sequence in $\wh G $ with limit set $L((\pi_k)_k) $. Let $O\in\g^*/G $ be the Kirillov orbit of some $\pi\in L((\pi_k)_k) $, $O_k$ the Kirillov orbit of $\pi_k$ for every $k$ and let $\ell \in O$. Then there exists for every $k $ an element $\ell_k\in O_k $, such that $\lim \limits_{k \to \infty} \ell_k=\ell $ in $\g^* $ (see \cite{leptin-ludwig}). One can assume that, passing to a subsequence if necessary, the sequence $(\g(\ell_k))_k $ converges in the subspace topology to a subalgebra $\u$ of  $ \g(\ell) $ and that there exists a number $d\in\N $, such that $\dim{O_k}=d $ for every $k\in\N $.
Then it follows from (\ref{desor}), that
\begin{eqnarray}\label{limorbit}
 L ((O_k)_k)=\lim \limits_{k \to \infty} \ell_k+\g(\ell_k)^\perp=\ell+\u^\perp\subset \g^*.
 \end{eqnarray}
Since $\g(\ell_k) $ contains $[\g,\g] $ for every $k $, the subspace $\u $ also contains $[\g,\g] $. Hence, the limit set $L((\pi_k)_k) $ in $\wh G $ of the sequence $(\pi_k)_k $ is the "affine" subset 
\begin{eqnarray*}
 L((\pi_k)_k)=\big\{ \big[\ch_q \otimes ind_P^G \chi_{\ell} \big]|~ q\in \u^\perp \big\}
 \end{eqnarray*}
for a polarization $\p$ in $\ell$ and $P:=exp(\p)$. \\
The observations above lead to the following proposition:
 
\begin{proposition}\label{threelim}
~\\
There are three different types of possible limit sets of the sequence $(O_k)_k$ of coadjoint orbits:
 \begin{enumerate}
\item The limit set $L((O_k)_k) $ is the singleton $O_\ell=\ell+\g(\ell)^{\perp} $, i.e. $\u=\g(\ell) $.
\item  The limit set $L((O_k)_k) $ is the  affine subspace $\ell+\u^\perp $ of characters of $\g $, i.e. $\big\langle{\ell},{[\g,\g]} \big\rangle=\{0\} $.
\item  The dimension of the orbit $O_\ell $ is strictly greater than $0 $ and strictly smaller than $d $. In this case
$$L((O_k)_k)=\bigcup_{\q\in\u^\perp} q+O_\ell,~~~\text{i.e.}~~~L((\pi_k)_k)=\bigcup_{q\in\u^\perp} \big[\ch_q\otimes ind_P^G \chi_{\ell}\big]$$
for a polarization $\p$ in $\ell$ and $P:=exp(\p)$. 
\end{enumerate}
 \end{proposition}

\subsection{The $C^*$-algebra $C^*(G/U,\ch_\ell) $}\label{CstUdef}
Let $\u\subset \g $ be an ideal of $\g $ containing $[\g,\g] $, $U:=exp(\u)$ and let $\ell\in\g^* $ such that $\langle{\ell},{[\g,\u]}\rangle=\{0\} $ and such that $\u\subset\g(\ell) $. Then the character $\ch_\ell $ of the group $U=\exp(\u) $ is $G $-invariant. One can thus define the involutive Banach algebra $\l1{G/U,\ch_\ell}$ as
\begin{eqnarray*}
 \l1{G/U,\ch_\ell}&:=& \bigg\{f:G\to \C|~ f\text{ measurable},~f(gu)=\ch_\ell(u\inv)f(g)~\forall~ g\in G\\
 & &~~~\forall~ u\in U,~\no f_1:=\int \limits_{G/U}\val{f(g)}~d\dot g<\iy \bigg\}.
 \end{eqnarray*}
The convolution
\begin{eqnarray*}
 f\ast f'(g):=\int \limits_{G/U}f(x)f'(x\inv g)~d\dot x~~~\forall~g \in G
 \end{eqnarray*}
and the involution
\begin{eqnarray*}
 f^*(g):=\ol{{f(g\inv)}}~~~\forall~g \in G
 \end{eqnarray*}
are well-defined for $f,f'\in \l1{G/U,\ch_\ell} $ and 
\begin{eqnarray*}
 \no{f\ast f'}_1\leq \|f\|_1 ~\|f' \|_1.
 \end{eqnarray*}
Moreover, the linear mapping
\begin{eqnarray*}
 p_{{G/U}}&:&L^1(G) \to \l1{G/U,\ch_\ell},~ \\
 &&p_{{G/U}}(F)(g):=\int \limits_U F(gu)\ch_\ell(u)~du
~~~\forall~ F\in \L1G~~\forall~g\in G
\end{eqnarray*}
is a surjective $* $-homorphism between the algebras $L^1(G) $ and $\l1{G/U,\ch_\ell} $. \\
Let 
\begin{eqnarray*}
 \wh G_{\u,\ell}~:=~\Big\{(\pi, \H_{\pi})\in \wh G|~\pi_{|_U}= {\ch_{\ell}}_{|_U} \Id_{\H_\pi} \Big\}.
 \end{eqnarray*}
Then $ \wh G_{\u,\ell} $ is a closed subset of $\wh G $, which can be identified with the spectrum of the algebra $\l1{G/U,\ch_\ell} $. Indeed it is  easy to see that every irreducible unitary  representation $(\pi,\H_\pi) \in \wh G_{\u,\ell}$ defines an irreducible representation $(\tilde \pi,\H_{\pi}) $ of the algebra $\l1{G/U,\ch_\ell} $ as follows:
\begin{eqnarray*}
 \tilde\pi \big(p_{G/U}(F) \big):=\pi(F)~~~\forall~ F\in L^1(G).
 \end{eqnarray*}
Similarly, if $(\tilde \pi,\H_{\tilde \pi}) $ is an irreducible unitary representation of $\l1{{G/U,\ch_\ell}} $ then 
\begin{eqnarray*}
 \pi:=\tilde\pi\ci p_{G/U}
 \end{eqnarray*}
defines an element of $\wh G_{\u,\ell} $. \\
Let $\s\subset \g $ be a subspace of $\g  $ such that $\g= \g(\ell)\oplus \s $. Since $\u $ contains $[\g,\g] $, it is easy to see that 
\begin{eqnarray}\label{whuell}
 \nn \wh G_{\u,\ell}= \big\{\big[\ch_q\otimes \pi_{\ell} \big]|~q\in(\u+\s)^\perp \big\},
 \end{eqnarray}
letting $\pi_{\ell}:= ind_P^G \chi_{\ell}$ for a polarization $\p$ in $\ell$ and $P:=exp(\p)$. \\
Denote by $C^*(G/U,\ch_\ell) $ the $C^* $-algebra of $L^1(G/U,\ch_\ell) $, whose spectrum can also be identified with $\wh G_{\u,\ell} $. \\
With $\pi_{\ell+q}= ind_P^G \chi_{\ell+q}$, the Fourier transform $\F $ defined by
\begin{eqnarray*}
 \F(a)(q):=\pi_{\ell+q}(a)~~~\forall~ q\in (\u+\s)^\perp
\end{eqnarray*} 
then maps the $C^* $-algebra $C^*(G/U,\ch_\ell) $ onto the algebra $C_0 \big((\u+\s)^\perp,\K(\H_{\pi_{\ell}}) \big) $ of the continuous mappings $\va:(\u+\s)^\perp\to  \K(\H_{\pi_{\ell}})$ vanishing at infinity with values in the algebra of compact operators on the Hilbert space of the representation $\pi_{\ell} $.\\

If one restricts $p_{G/U} $ to the Fr\'echet algebra $\S(G) \subset L^1(G)$, its image will be the Fr\'echet algebra 
\begin{eqnarray*}
\S(G/U,\ch_\ell)&=&\big\{f\in \l1{G/U,\ch_\ell}|~f~\text{smooth and for every subspace}~\s' \subset \g~\text{with}\\
&&~~\g=\s'\oplus \u~\text{and for}~S'=exp(\s'),~ f_{|_{S'}} \in \S(S') \big\}.
\end{eqnarray*}
 
\section{Conditions 1, 2 and 3(a)}\label{cond 1,2,3a}
Now, to start with the proof of the above listed conditions, the families of sets $(S_i)_{i \in \{0,...,r\}}$ and $(\Gamma_i)_{i \in \{0,...,r\}}$ are going to be defined and the Properties 1, 2 and 3(a) are going to be checked.

In order to be able to do this, one needs to construct a polarization $\p_{\ell}^V$ for $\ell \in \g^*$ as follows: \\
~\\
Fix once and for all a Jordan-H\"older basis $\{H_1,...,H_n\}$ of $\g$, in such a way that $\g_i:= span\{H_i,...,H_n\}$ for $i \in \{0,...,n\}$ is an ideal in $\g $. Since $\g $ is two-step nilpotent, one can first choose a basis $\{H_{\tilde{n}},\cdots, H_n\} $ of $[\g,\g]$ and then add the vectors $H_1,\cdots, H_{\tilde{n}-1} $ to obtain a basis of $\g$. Let 
$$I_{\ell}^{Puk}:=\{ i \leq n |~\g(\ell) \cap \g_i = \g(\ell) \cap \g_{i+1} \}$$ be the Pukanszky index set for $\ell \in \g^*$. The number of elements $\big|I^{Puk}_\ell \big| $ of $ I^{Puk}_\ell$ is the dimension of the orbit $O_\ell $ of $\ell $. \\
Moreover, if one denotes by $\g_i(\ell_{| \g_i})$ the stabilizer of $\ell_{| \g_i}$ in $\g_i$,
$$\p_{\ell}^V:= \sum \limits_{i=1}^{n} \g_i(\ell_{| \g_i})$$
is the Vergne polarization of $\ell$ in $\g$. Its construction will now be analyzed by a method developed in \cite{ludwig-zahir}, Section 1: \\
Let $\ell \in \g^*$. Then choose the greatest index $j_1(\ell) \in \{1,...,n\}$ such that \mbox{$H_{j_1(\ell)} \not \in \g(\ell)$} and let $Y_1^{V,\ell}:=H_{j_1(\ell)}$. Furthermore, choose the index $k_1(\ell) \in \{1,...,n\}$ such that $\langle \ell,[H_{k_1(\ell)},H_{j_1(\ell)}]\rangle \not=0$ and $\langle \ell, [H_i,H_{j_1(\ell)}]\rangle=0$ for all $i>k_1(\ell)$ and let $X_1^{V,\ell}:=H_{k_1(\ell)}$. \\
Next, let $\g^{1,\ell}:=\big\{U \in \g |~\big\langle \ell,\big[U,Y_1^{V,\ell}\big]\big\rangle=0\big\}$. Then $\g^{1,\ell}$ is an ideal in $\g$ which does not contain $X_1^{V,\ell}$, and $\g=\R X_1^{V,\ell} \oplus \g^{1,\ell}$. Now, the Jordan-H\"older basis will be changed, taking out $H_{k_1(\ell)}$: \\
Consider the Jordan-H\"older basis $\{H_1^{1,\ell},...,H_{k_1(\ell)-1}^{1,\ell},H_{k_1(\ell)+1}^{1,\ell},...,H_n^{1,\ell}\}$ of $\g^{1,\ell}$ with 
$$H_i^{1,\ell}:=H_i~~\forall~i >k_1(\ell)~~~\text{and}~~~
H_i^{1,\ell}:=H_i- \frac{\big\langle \ell,\big[H_i,Y_1^{V,\ell}\big]\big\rangle X_1^{V,\ell}}{\big\langle \ell,\big[X_1^{V,\ell},Y_1^{V,\ell}\big] \big\rangle}~~\forall~i<k_1(\ell).$$
Then, choose the greatest index $j_2 (\ell)\in \{1,...,k_1(\ell)-1,k_1(\ell)+1,...,n\}$ in such a way that \mbox{$H_{j_2(\ell)}^{1,\ell} \not \in \g^{1,\ell}(\ell_{|\g^{1,\ell}})$} and define $Y_2^{V,\ell}:=H_{j_2(\ell)}^{1,\ell}$. Like above, choose \mbox{$k_2(\ell) \in \{1,...,k_1(\ell)-1,k_1(\ell)+1,...,n\}$} such that $\langle \ell,[H_{k_2(\ell)}^{1,\ell},H_{j_2(\ell)}^{1,\ell}]\rangle \not=0$ and that $\langle \ell, [H_i^{1,\ell},H_{j_2(\ell)}^{1,\ell}] \rangle=0$ for all $i>k_2(\ell)$ and set $X_2^{V,\ell}:=H_{k_2(\ell)}^{1,\ell}$. \\
Iterating this procedure, one gets sets $\big\{Y_1^{V,\ell},...,Y_d^{V,\ell}\big\}$ and $\big\{X_1^{V,\ell},...,X_d^{V,\ell}\big\}$ for $d \in \{0,...,[\frac{n}{2}]\}$ with the properties $$\p_{\ell}^V=span\big\{Y_1^{V,\ell},...,Y_d^{V,\ell}\big\} \oplus \g(\ell)$$ and $$\big\langle \ell ,\big[X_i^{V,\ell},Y_i^{V,\ell}\big]\big\rangle \not=0,~~~\big\langle \ell,\big[X_i^{V,\ell},Y_j^{V,\ell}\big] \big\rangle =0 ~~~ \forall ~i \not= j \in \{1,...,d\}~~~\text{and}$$ 
$$\big\langle \ell ,\big[Y_i^{V,\ell},Y_j^{V,\ell}\big]\big\rangle =0~~~ \forall~ i,j \in \{1,...,d\}.$$

Now, let \begin{eqnarray}\label{inddouble}
 \nn J(\ell) ~:= ~\{j_1(\ell),\cdots,j_d(\ell) \}~~~\text{and}~~~K(\ell):=\{k_1(\ell),\cdots, k_d(\ell)\}. 
\end{eqnarray}
Then
\begin{eqnarray}\label{pukdecom}
 \nn I^{Puk}_\ell ~= ~J(\ell)~\dot\cup~ K(\ell)~~~\text{and}~~~j_1(\ell)>\cdots>j_d(\ell).
\end{eqnarray}
It is easy to see that the index sets $I^{Puk}_\ell, J(\ell) $ and $K(\ell) $ are the same on every coadjoint orbit (see \cite{ludwig-zahir}) and can therefore also be denoted by $I^{Puk}_O,~ J(O) $ and $K(O)$ if $\ell$ is located in the coadjoint orbit $O$.

Now, for the parametrization of $\g^*/G$ and thus of $\wh{G}$ and for the choice of the in Property 1(b) required concrete realization of a representation, let $O \in \g^*/G$. A theorem of L.Pukanszky (see \cite{puk}, Part II, Chapter I.3 or \cite{ludwig-zahir}, Corollary 1.2.5) states that there exists one unique $\ell_O \in O$ such that $\ell_O(H_i)=0$ for every index $i \in I_O^{Puk}$. So, choose this $\ell_O$, let $P_{\ell_O}^V:=\exp(\p_{\ell_O}^V)$ and define the irreducible unitary representation 
$$\pi_{\ell_O}^V:=ind_{P_{\ell_O}^V}^G \chi_{\ell_O}$$ 
associated to the orbit $O$ and acting on $L^2\big(G/P_{\ell_O}^V,\chi_{\ell_O}\big)\cong L^2(\R^d)$.\\
~\\
\newpage
Next, one has to construct the demanded sets $\Gamma_i$ for $i \in \{0,...,r\}$: \\
For this, define for a pair of sets $(J,K)$ such that $J,K\subset \{1,\cdots, n\}, ~\val J=\val K$ and $J\cap K=\es $
the subset $ (\g^*/G)_{(J,K)} $ of $\g^*/G $ by 
\begin{eqnarray*}
 (\g^*/G)_{(J,K)} :=\{O\in \g^*/G|~ (J,K)=(J(O), K(O))\}.
 \end{eqnarray*}
Moreover, let 
 $$\M:=\{(J,K)|~J,K \subset \{1,...,n\},~J \cap K = \emptyset,~|J|=|K|,~(\g^*/G)_{(J,K)} \not= \emptyset \}$$
 and
$$(\g^*/G)_{2d}:=\big\{O \in \g^*/G|~\big|I_O^{Puk}\big|=2d\big\}.$$ 

Then $$(\g^*/G)_{2d}~=~\bigcup_{\substack{(J,K):~J,K \subset \{1,...,n\},\\|J|=|K|=d,~J\cap K=\emptyset}} (\g^*/G)_{(J,K)}$$
and $$\g^*/G~=~\bigcup \limits_{d \in \{0,...,[\frac{n}{2}]\}} (\g^*/G)_{2d}
~=~\bigcup_{(J,K)\in\M}(\g^*/G)_{(J,K)}.$$
Now, an order on the set $\M$ shall be introduced.\\
First, if $|J|=|K|=d$, $|J'|=|K'|=d'$ and $d<d'$, then the pair $(J,K)$ is defined to be smaller than the pair $(J',K')$: $(J,K)<(J',K')$.\\
If $|J|=|K|=|J'|=|K'|=d$, $J=\{j_1,...,j_d\}$, $J'=\{j'_1,...,j'_d\}$ and $j_1<j'_1$, the pair $(J,K)$ is again defined to be smaller than $(J',K')$. \\
Otherwise, if $j_1=j'_1$, one has to consider $K=\{k_1,...,k_d\}$ and $K'=\{k'_1,...,k'_d\}$ and here again, compare the first elements $k_1$ and $k'_1$:
So, if $j_1=j'_1$ and $k_1<k'_1$, again $(J,K)<(J',K')$. \\
But if $k_1=k'_1$, one compares $j_2$ and $j'_2$ and continues in that way. \\
If $r+1=\val\M $, one can identify the ordered set $\M $ with the interval $\{0,\cdots, r\} $ and assign to each such pair $(J,K)\in \M$ a number $i_{JK} \in \{0,...,r\}$. \\
Finally, one can therefore define the sets $\Gamma_{i_{JK}}$ and $S_{i_{JK}}$ as
$$\Gamma_{i_{JK}}~:=~\big\{\big[\pi_{\ell_O}^V\big]|~O \in (\g^*/G)_{(J,K)}\big\}~~~\text{and}$$
$$S_{i_{JK}}~:=~\bigcup \limits_{i \in \{0,...,i_{JK}\}} \Gamma_{i}.$$
Then obviously, the family $(S_i)_{i \in \{0,...,r\}}$ is an increasing family in $\wh{G}$. \\
Furthermore, the set $S_i$ is closed for every $i \in \{0,...,r\}$. This can easily be deduced from the definition of the index sets $J(\ell)$ and $K(\ell) $. The  indices $j_m(\ell)$ and $k_m(\ell)$ for $m \in \{1,...,d\}$ are chosen in such a way that they are the largest to fulfill a condition of the type $\langle \ell, [H_{j_m(\ell)}^{m-1,\ell},\cdot]\rangle \not=0$ or $\langle \ell, [H_{k_m(\ell)}^{m-1,\ell}, \cdot] \rangle \not= 0$, respectively.\\
In addition, the sets $\Gamma_i$ are Hausdorff. For this, let $i=i_{JK}$ for $(J,K) \in \M$ and $(O_k)_k$ in $(\g^*/G)_{(J,K)}$ a sequence of orbits such that the sequence $\big(\big[\pi_{\ell_{O_k}}^V\big]\big)_k$ converges in $\Gamma_i$, i.e. $(O_k)_k$ converges in $(\g^*/G)_{(J,K)}$ and thus has a limit point $O$ in $(\g^*/G)_{(J,K)}$. If now $O_k \ni \ell_k \overset{k \to \infty}{\longrightarrow} \ell \in O$, then by (\ref{limorbit}), it follows that the limit $\u $ of the sequence $(\g(\ell_k))_k$ is equal to $\g(\ell)$. Therefore, the sequence $(O_k)_k$ and thus also the sequence $\big(\big[\pi_{\ell_{O_k}}^V\big]\big)_k$ have unique limits and hence $\Gamma_i$ is Hausdorff.\\
Moreover, one can still observe that for $d=0$ the choice $J=K=\emptyset$ represents the only possibility to get $|J|=|K|=d$. So, the pair $(\emptyset, \emptyset)$ is the first element in the above defined order and therefore corresponds to $0$. Thus 
$$\Gamma_0~=~\big\{\big[\pi_{\ell_O}^V\big]|~I_O^{Puk}=\emptyset\big\},$$
which is equivalent to the fact that $\g(\ell_O)=\g$ which again is equivalent to the fact that every $\pi_{\ell_O}^V \in \Gamma_0$ is a character. Hence, $S_0=\Gamma_0$ is the set of all characters on $\g$, as demanded. \\
~\\
Since one can identify the quotient space $G/P_{\ell_O}^V$ with $\R^d$ by means of the subspace $\s_{\ell_O}=\span{X_1^{V,\ell_O},...,X_d^{V,\ell_O}}$, one can identify the Hilbert space $L^2\big(G/P_{\ell_O}^V,\chi_{\ell_O}\big)$ with $L^2(\R^d) $ as in (\ref{indident}) and one can suppose that the representation $\pi_{\ell_O}^V$ acts on the Hilbert space $L^2(\R^d)$ for every $O \in (\g^*/G)_{2d}$. \\
~\\
Hence, the first condition is fulfilled. For the proof of the Properties 2 and 3(a), a proposition will be shown: 

\begin{proposition}
~\\
For every $a \in C^*(G)$ and every $(J,K) \in \M$ with $|J|=|K|=d \in \{0,...,[\frac{n}{2}]\}$, the mapping 
$$\Gamma_{i_{JK}} \to L^2(\R^d), ~\gamma \mapsto \F(a)(\gamma)$$ 
is norm continuous and the operator $\F(a)(\gamma)$ is compact for all $\gamma \in \Gamma_{i_{JK}}$.
\end{proposition}

Proof: \\
The compactness follows directly from a general theorem which can be found in \cite{cor-green} (Chapter 4.2) or \cite{puk} (Part II, Chapter II.5) and states that the $C^*$-algebra $C^*(G)$ of every connected nilpotent Lie group $G$ fulfills the CCR condition, i.e. the image of every irreducible representation of $C^*(G)$ is a compact operator. \\
~\\
Next, let $d \in \{0,...,[\frac{n}{2}]\}$ and $(J,K) \in \M$ such that $|J|=|K|=d$. \\
First, one has to observe that the polarization $\p_{\ell}^V$ is continuous in $\ell$ on the set $\big\{\ell_{O'}|~ O' \in (\g^*/G)_{(J,K)} \big\}$. This can be seen by the construction of the vectors $\big\{Y_1^{V,\ell},...,Y_d^{V,\ell} \big\}$. \\
Now, let $(O_k)_k$ be a sequence in $(\g^*/G)_{(J,K)}$ and $O \in (\g^*/G)_{(J,K)}$ such that $\big[\pi_{\ell_{O_k}}^V\big] \overset{k \to \infty}{\longrightarrow} \big[\pi_{\ell_O}^V\big]$ and let $a \in C^*(G)$. Then $\ell_{O_k} \overset{k \to \infty}{\longrightarrow} \ell_O$ and by the observation above, the associated sequence of polarizations $\big(\p_{\ell_{O_k}}^V\big)_k$ converges to the polarization $\p_{\ell_O}^V$. By Theorem 2.3 in \cite{hedi}, thus $\pi_{\ell_{O_k}}^V(a) \overset{k \to \infty}{\longrightarrow} \pi_{\ell_O}^V(a)$ in the operator norm. \\
\qed
~\\
~\\
Since $C^*(G)$ is obviously separable, this proposition proves the desired Properties 2 and 3(a) and hence, it remains to show Property 3(b): 
\section{Condition 3(b)}\label{Step 2}
\subsection{Introduction to the setting}\label{intro}
For simplicity, in the following, the representations will be identified with their equivalence classes. \\
~\\
Let $d \in \{0,...,[\frac{n}{2}]\}$ and $(J,K) \in \M$ with $|J|=|K|=d$. Furthermore, fix $i=i_{JK} \in \{0,...,r\}$. \\
~\\
Let $\big(\pi_k^V\big)_k=\big(\pi_{\ell_{O_k}}^V\big)_k$ be a sequence in $\Gamma_i$ whose limit set is located outside $\Gamma_i$. Since \mbox{every} converging sequence has a properly converging subsequence, it will be assumed that $\big(\pi_k^V\big)_k$ is properly converging and the transition to a subsequence will be omitted. \\ 
The corresponding sequence of coadjoint orbits $(O_k)_k$ is contained in $(\g^*/G)_{(J,K)}$ and in particular every $O_k$ has the same dimension $2d$. Moreover, it converges properly to a set of orbits $L((O_k)_k)$.\\
In addition, since $S_i$ is closed, the limit set $L \big( \big(\pi_k^V \big)_k\big)$ of the sequence $\big(\pi_k^V\big)_k$ is contained in $S_{i-1}=\bigcup \limits_{j \in \{0,...,i-1\}} \Gamma_j$ and therefore for every element $O \in L((O_k)_k)$ there exists a pair $(J_O,K_O)< (J,K)$ such that $\pi_{\ell_O}^V \in \Gamma_{i_{J_O K_O}}$ or equivalently, $O \in (\g^*/G)_{(J_O,K_O)}$.

\subsection{Changing the Jordan-H\"older basis.}
Let $\tilde{\ell} \in \tilde{O} \in L((O_k)_k)$. Then, there exists a sequence $\big(\tilde{\ell}_k \big)_k$ in $O_k$ such that $\tilde{\ell}= \lim \limits_{k \to \infty} \tilde{\ell}_k$. \\
Since one is interested in the orbits $O_k=\tilde{\ell_k}+\g \big(\tilde{\ell}_k \big)^{\perp}$, one can change the sequence $\big(\tilde{\ell}_k \big)_k$ to a sequence $(\ell_k)_k$ by letting $\ell_k(A)=0$ for every $A \in \g \big(\tilde{\ell}_k \big)^{\perp}=\g(\ell_k)^{\perp}$. \\
Thus, one obtains another converging sequence $(\ell_k)_k$ in $(O_k)_k$ whose limit $\ell$ is located in an orbit $O \in L((O_k)_k)$. \\
As above, one can suppose that the subalgebras $(\g(\ell_k))_k $ converge to a subalgebra $ \u $, whose corres\-ponding Lie group $exp (\u)$ is denoted by $U$. These subalgebras $ \g(\ell_k) $ can be written as
$$\g(\ell_k)~=~[\g,\g] \oplus\s_k,$$
where $ \s_k\subset [\g,\g] ^\perp $. In addition, let $ \n_{k,0}$  be the kernel of $ \ell_k{\res{[\g,\g] }} $ and $ \s_{k,0} $ the kernel of $ \ell_k{\res{\s_k}}$ for all $k \in \N$. One can assume  that $ \s_{k,0}\ne\s_k $ and choose $ T_k \in \s_k $ orthogonal to $ \s_{k,0} $ of length 1. The case   $ \s_{k,0}=\s_k$ for  $k\in\N$, being easier, will be omitted. \\
Similarly, choose $ Z_k\in [\g,\g]  $ orthogonal to $ \n_{k,0} $ of length 1.
One sees that such a $Z_k$ must exist: If $ \ell_k{\res{[\g,\g] }}=0$ for $k\in\N$, then $\pi_{\ell_{O_k}}^V$ is a character and thus contained in $S_0=\Gamma_0$.  But $S_0=\Gamma_0$ is closed and thus $\big(\pi_{\ell_{O_k}}^V\big)_k$ cannot have a limit set outside $\Gamma_0$. \\
~\\
Furthermore, let $ \r_k=\g(\ell_k)^\perp\subset\g $. \\
One can assume that, passing to a subsequence if necessary, $\lim \limits_{k \to \infty} Z_k=:Z$, $\lim \limits_{k \to \infty} T_k =:T $ and $ \lim \limits_{k \to \infty} \r_k=:\r $ exist. \\
~\\
Now, new polarizations $\p_k$ in $\ell_k$ are  needed: \\
The restriction to $ \r_k $ of the skew-form $ B_k:=B_{\ell_k} $ defined in Chapter \ref{pre} is non-degenerate on $ \r_k $ and there exists an invertible endomorphism $ S_k $ of $ \r_k $ such that
$$\langle{x},{S_k(x')}\rangle ~= ~B_k(x,x')~~~\forall~ x,x'\in \r_k.$$
Then $ S_k $ is skew-symmetric, i.e. $ S_k^t=-S_k $, and with the help of Lemma \ref{lemma} one can decompose $\r_k $ into an orthogonal direct sum
$$\r_k ~=~\sum_{j=1}^d V_k^j$$
of two-dimensional $ S_k $-invariant subspaces.
Choose an orthonormal basis $ \{X_j^k, Y_j^k\} $ of $ V^j_k $. 
Then,
\begin{eqnarray*}
&&[X_i^k,X^k_j] \in \n_{k,0}~~~\forall~ i, j \in \{1,...,d\}, \\
&&[Y_i^k,Y_j^k] \in \n_{k,0}~~~\forall~ i, j \in \{1,...,d\}~~~\text{and}\\
&&[X^k_i,Y_j^k] ~=~\de_{i,j}c_j^k Z_k\text{ mod }\n_{k,0}~~~\forall~ i,j \in \{1,...,d\},
\end{eqnarray*}
whereat $0 \ne c^k_j\in\R $ and $\sup \limits_{k \in \N} c^k_j < \infty$ for every $j \in \{1,...,d\}$. \\
Again, by passing to a subsequence if necessary, the sequence $(c^k_j)_k$ converges for every $j \in \{1,...,d\}$ to some $c_j$. \\
Since $X_j^k, Y_j^k \in \r_k$ and $\ell_k(A)=0$ for every $A \in \r_k$, $\ell_k(X_j^k)=\ell_k(Y_j^k)=0$ for all $j \in \{1,...,d\}$. Furthermore, one can suppose that the sequences $ (X^k_j)_k,( Y^k_j)_k $ converge in $ \g $ to vectors $ X_j, Y_j $ which form a basis modulo $ \u $ in $ \g $.\\
It follows that
\begin{eqnarray*}
 \langle{\ell_k},{[X^k_j,Y_j^k] }\rangle &= &c^k_j\la_k,~\text{where}\\
 \la_k&=&\langle{\ell_k},{Z_k}\rangle \overset{k \to \infty}{\longrightarrow} \langle{\ell},{Z}\rangle~=: \lambda.
\end{eqnarray*}
As $Z_k$ was chosen orthogonal to $\n_{k,0}$, $\la_k \not=0$ for every $k$. \\

Now, let 
$$\p_k ~:= ~\text{span}\{Y_1^k,\cdots, Y^k_d,\g(\ell_k)\}$$
and $P_k:=\exp(\p_k)$. Then $ \p_k $ is a polarization at $ \ell_k $. Furthermore, define the re\-presentation $\pi_k$ as
$$\pi_k ~:= ~\ind_{P_k}^G\ch_{\ell_k}.$$
Then, since $\pi_k$, as well as $\pi_k^V$ are induced representations of polarizations and of the characters $\chi_{\ell_k}$ and $\chi_{\ell_{O_k}}$, whereat $\ell_k$ and $\ell_{O_k}$ lie in the same coadjoint orbit $O_k$, the two representations are equivalent, as observed in Chapter \ref{pre}.  \\
~\\
Let $ \a_k:=\n_{k,0}+\s_{k,0} $. Then $ \a_k $ is an ideal of $ \g $ on which  $ \ell_k $ is 0. Therefore, the normal subgroup $\exp(\a_k)$ is contained in the kernel of the representation $ \pi_{k} $.
Moreover, let $ \a:=\lim \limits_{k \to \infty} \a_k$. \\
In addition, let $p \in \N$, $\tilde{p} \in \{1,...,p\}$ and let $\{A^k_1,\cdots,A^k_{\tilde{p}} \}$ denote an orthonormal basis of $\n_{k,0}$, the part of $\a_k$ which lies inside $[\g,\g]$, and $\{A^k_{\tilde{p}+1},\cdots,A^k_p\}$ an orthonormal basis of $\s_{k,0}$, the part of $\a_k$ outside $[\g,\g]$. Then $\{A^k_1,\cdots,A^k_p\}$ is an orthonormal basis of $ \a_k$ and as above, one can assume that $ \lim \limits_{k \to \infty} A^k_j=A_j $ exists for all $j \in \{1,...,p\}$. \\
Now, for every $k \in \N$ one can take as an orthonormal basis for $\g$ the set of vectors
$$\{X^k_1,\cdots, X^k_d,Y^k_1,\cdots, Y^k_d,T_k,Z_k,A_1^k,...,A_p^k\}$$
as well as the set 
$$\{X_1,\cdots, X_d,Y_1,\cdots, Y_d,T,Z,A_1,...,A_p\}. $$
This gives the following Lie brackets:
\begin{eqnarray*}
 &[X^k_i,Y^k_j]&  = \de_{i,j}c^k_j Z_k\text{ mod }\a_k,\\
&[X^k_i,X^k_j]& =0 \text{ mod }\a_k~~~\text{and}\\
&[Y^k_i,Y^k_j]& =0 \text{ mod }\a_k.
\end{eqnarray*}
The vectors $ Z_k $ and $ T_k $ are central modulo $ \a_k $. \\
~\\
Before starting the analysis of $(\pi_k)_{k \in \N}$, some notations have to be introduced: 

\subsection{Definitions}\label{definitions}
Choose for $j \in \{1,...,d\}$ the Schwartz functions $\eta_j \in \cal{S}(\R)$ such that $\| \eta_j \|_{L^2(\R)}=1$ and $\| \eta_j \|_{L^{\infty}(\R)}\klgl 1.$ \\
Furthermore, for $x_1,...,x_d,y_1,...,y_d,t,z,a_1,\cdots,a_p \in \R$, write
$$(x)_k:=(x_1,...,x_d)_k:=\sum_{j=1}^d x_j X^k_j,~~(y)_k~:=(y_1,...,y_d)_k:=\sum_{j=1}^d y_j Y^k_j,~~(t)_k:=tT_k,$$
$$(z)_k:=z Z_k,~~(\dot{a})_k:=(a_1,...,a_{\tilde{p}})_k:=\sum_{j=1}^{\tilde{p}} a_j A^k_j,~~(\ddot{a})_k:=(a_{\tilde{p}+1},...,a_p)_k:=\sum_{j=\tilde{p}+1}^p a_j A^k_j$$
$$\text{and}~~(a)_k~:=(\dot{a},\ddot{a})_k=(a_1,...,a_p)_k=\sum_{j=1}^p a_j A^k_j,$$
whereat $(\cdot,...,\cdot)_k$ is defined to be the $d$-, $\tilde{p}$-, $(p-\tilde{p})$- or the $p$-tuple with respect to the bases $\{X_1^k,...,X_d^k\}$, $\{Y_1^k,...,Y_d^k\}$, $\{A_1^k,...,A_{\tilde{p}}^k\}$, $\{A_{\tilde{p}+1}^k,...,A_p^k\}$ and $\{A_1^k,...,A_p^k\}$, respectively, and let
\begin{eqnarray*}
(g)_k&:=&(x_1,...,x_d,y_1,...,y_d,t,z,a_1,...,a_p)_k~:=~((x)_k,(y)_k,(t)_k,(z)_k,(\dot{a})_k,(\ddot{a})_k)\\
&=&((x)_k,(h)_k)\\
&=&\sum \limits_{j=1}^d x_j X^k_j~+~ \sum
\limits_{j=1}^d y_j Y_j^k~+~t T_k~+~z Z_k~+~\sum \limits_{j=1}^p a_j A_j^k,
\end{eqnarray*}
where $(h)_k$ is in the polarization $\p_k$ and the $(2d+2+p)$-tuple $(\cdot,...,\cdot)_k$ is regarded with respect to the basis $\{X_1^k,...,X_d^k,Y_1^k,...,Y_d^k,T_k,Z_k,A_1^k,...,A_p^k\}$. \\
Moreover, define the limits
$$(x)_{\infty}:=(x_1,...,x_d)_{\infty}:= \lim \limits_{k \to \infty} (x)_k=\sum \limits_{j=1}^d x_j X_j,
~~(y)_{\infty}:=(y_1,...,y_d)_{\infty}:= \lim \limits_{k \to \infty} (y)_k=\sum \limits_{j=1}^d y_j Y_j,$$
$$(t)_{\infty}:=\lim \limits_{k \to \infty} (t)_k=t T,~~(z)_{\infty}:=\lim \limits_{k \to \infty} (z)_k=z Z, ~~(\dot{a})_{\infty}:=(a_1,...,a_{\tilde{p}})_{\infty}:= \lim \limits_{k \to \infty} (\dot{a})_k=\sum \limits_{j=1}^{\tilde{p}} a_j A_j,$$
$$(\ddot{a})_{\infty}:=(a_{\tilde{p}+1},...,a_p)_{\infty}:= \lim \limits_{k \to \infty} (\ddot{a})_k=\sum \limits_{j=\tilde{p}+1}^{p} a_j A_j,$$
$$(a)_{\infty}:=(\dot{a},\ddot{a})_{\infty}=(a_1,...,a_p)_{\infty}:= \lim \limits_{k \to \infty} (a)_k=\sum \limits_{j=1}^p a_j A_j~~~\text{and}$$
$$(g)_{\infty}:=(x,y,t,z,\dot{a},\ddot{a})_{\infty}~:=~ \lim \limits_{k \to \infty} (g)_k~=~\sum \limits_{j=1}^d x_j X_j~+~ \sum
\limits_{j=1}^d y_j Y_j~+~t T~+~z Z~+~\sum \limits_{j=1}^p a_j A_j.$$
~\\
Now, the representations $(\pi_k)_{k \in \N}$ can be computed: 

\subsection{Formula for $\pi_k$}\label{Kapitel Pik-Formel}
Let $f \in L^1(G)$. \\
With $\rho_k:= \langle \ell_k, T_k \rangle$, $c^k:=(c_1^k,...,c_d^k)$ and $(s)_k:=(s_1,...,s_d)_k=\sum \limits_{j=1}^d s_j X_j^k$ for $s_1,...,s_d \in \R$, where again $(\cdot,...,\cdot)_k$ is the $d$-tuple with respect to the basis $\{X_1^k,...,X_d^k\}$, as in (\ref{exppi}), the representation $ \pi_k $ acts on $L^2(G/P_k, \chi_{\ell_k})$ in the following way:
\begin{eqnarray}\label{replk}
\nn \pi_k\big((g)_k\big)\xi\big((s)_k\big) &=&\xi\big((g)_k\inv\cdot (s)_k\big)\\
\nn&=&e^{2 \pi i\langle{\ell_k},{-(y)_k-(t)_k-(z)_k-(\dot{a})_k-(\ddot{a})_k+[(s)_k-\frac{1}{2}(x)_k,(y)_k]-\frac{1}{2}[(x)_k,(s)_k]}\rangle}\xi\big((s-x)_k\big)\\
\nn&=&e^{2 \pi i\big(- t\rho_k -z\la_k+\sum \limits_{j=1}^d \la_k c^k_j(s_j-\frac{1}{2}x_j)y_j\big)}
\xi\big((s-x)_k\big) \\
&=&e^{2\pi i (-t \rho_k-z\la_k+ \la_k c^k((s)_k-\frac{1}{2}(x)_k)(y)_k)} \xi\big((s-x)_k\big),
\end{eqnarray}
since $\ell_k(Y_j^k)=0$ for all $j \in \{1,...,d\}$. \\

Now, identify $G$ with $\R^d \times \R^d \times \R \times \R \times \R^{\tilde{p}} \times \R^{p-{\tilde{p}}} \cong \R^{2d+2+p}$, let $\xi \in L^2(\R^d)$ and $s \in \R^d$. Moreover, identify $\pi_k$ with a representation acting on $L^2(\R^d)$ which will also be called $\pi_k$. To stress the dependence on $k$ of the above fixed function $f \in L^1(G)$, denote by $f_k \in L^1(\R^{2d+2+p})$ the function $f$ applied to an element in the $k$-basis: $$f_k(g):=f((g)_k).$$ Then,
\begin{eqnarray}\label{Pik-Formel1}
\nn &&\pi_k(f)\xi(s)\\
\nn~\\
\nn&=& \int \limits_{\R^d \times \R^d \times \R \times \R \times \R^{\tilde{p}}\times \R^{p-\tilde{p}}} f_k(g) \pi_k(g) \xi(s)~dg \\
\nn~\\
\nn&=& \int \limits_{\R^d \times \R^d \times \R \times \R \times \R^{\tilde{p}}\times \R^{p-\tilde{p}}}
f_k(x,y,t,z,\dot{a},\ddot{a}) e^{2\pi i (-t \rho_k-z\la_k+ \la_k c^k(s-\frac{1}{2}x)y)} \xi(s-x)~d(x,y,t,z,\dot{a},\ddot{a})\\
\nn~\\
\nn&=& \int \limits_{\R^d \times \R^d \times \R \times \R \times \R^{\tilde{p}}\times \R^{p-\tilde{p}}}
f_k(s-x,y,t,z,\dot{a},\ddot{a}) e^{2\pi i (-t \rho_k-z\la_k+
\frac{1}{2} \la_k c^k(s+x)y)} \xi(x)~d(x,y,t,z,\dot{a},\ddot{a}) \\
\nn~\\
&=&  \int \limits_{\R^d} \hat f_k^{2,3,4,5,6}\Big(s-x,-\frac{\lambda_k c^k}{2}(s+x),\rho_k,\lambda_k,0,0\Big)
 \xi(x)~dx,
\end{eqnarray}

where $\hat{f}_k^{2,3,4,5,6}$ denotes the Fourier transform in the 2nd, 3rd, 4th, 5th and 6th variable.

\subsection{First case}\label{first}
First consider the case that $L((O_k)_k)$ consists of one single limit point $O$. \\
In this case, for every $k$, 
$$2d~=~dim(O_k)~=~dim(O).$$
Thus, the regarded situation occurs if and only if $\lambda \not=0$ and $c_j \not= 0$ for every $j \in \{1,...,d\}$. \\

Consider again the above chosen sequence $(\ell_k)_k$ which converges to $\ell \in O$. As the dimensions of the orbits $O_k$ and $O$ are the same, there exists a subsequence of $(\ell_k)_k$ (which will also be denoted by $(\ell_k)_k$ for simplicity) such that $\p:= \lim \limits_{k \to \infty} \p_{\ell_k}^V$ is a polarization for $\ell$, but not necessarily the Vergne polarization. Moreover, define $P:=exp(\p)= \lim \limits_{k \to \infty}P_{\ell_k}^V$ and let
$$\pi ~:= ~\ind_{P}^G\ch_{\ell}.$$
Now, if one identifies the Hilbert spaces $\H_{\pi_{\ell_k}^V}$ and $\H_{\pi}$ of $\pi_{\ell_k}^V=\ind_{P_{\ell_k}^V} ^G\ch_{\ell_k}$ and $\pi$ with $L^2(\R^d)$, from \cite{hedi}, Theorem 2.3, one can conclude that 
$$\big\| \pi_{\ell_k}^V(a)-\pi(a) \big\|_{op}~=~\Big\|\ind_{P_{\ell_k}^V} ^G\ch_{\ell_k}(a)-\ind_{P}^G\ch_{\ell}(a)\Big\|_{op} \overset{k \to \infty}{\longrightarrow}0~~~\forall~ a \in C^*(G).$$ 
Since $\pi$ and $\pi_{\ell}^V=\ind_{P_{\ell}^V} ^G\ch_{\ell}$ are both induced representations of polarizations and of the same character $\chi_{\ell}$, they are equivalent and hence, there exists a unitary intertwining operator 
$$F: \H_{\pi_{\ell}^V} \cong L^2(\R^d) \to \H_{\pi} \cong L^2(\R^d)~~\text{such that}~~F \circ \pi_{\ell}^V(a) = \pi(a) \circ F~~~\forall ~a \in C^*(G).$$
Moreover, the two representations $\pi_k^V=\pi_{\ell_{O_k}}^V=ind_{P_{\ell_{O_k}}^V}^G \chi_{\ell_{O_k}}$ and $\pi_{\ell_k}^V=ind_{P_{\ell_k}^V}^G \chi_{\ell_k}$ are equivalent for every $k \in \N$ because $\ell_{O_k}$ and $\ell_k$ are located in the same coadjoint orbit $O_k$ and $\p_{\ell_{O_k}}^V$ and $\p_{\ell_k}^V$ are polarizations. Thus there exist further unitary intertwining operators
$$F_k: \H_{\pi_k^V} \cong L^2(\R^d) \to \H_{\pi_{\ell_k}^V} \cong L^2(\R^d)~~\text{with}~~F_k \circ \pi_k^V(a) = \pi_{\ell_k}^V(a) \circ F_k~~~\forall~ a \in C^*(G).$$
Now, define the required operators $\tilde{\nu}_k$ as 
$$\tilde{\nu}_k(\varphi):=F_k^* \circ F \circ \varphi \big(\pi_{\ell}^V \big) \circ F^* \circ F_k ~~~\forall ~\varphi \in CB(S_{i-1}),$$ 
which makes sense since $\pi_{\ell}^V$ is a limit point of the sequence $\big(\pi_k^V \big)_k$ and hence contained in $S_{i-1}$, as seen in Chapter \ref{intro}. \\
As $\varphi \big(\pi_{\ell}^V \big) \in \B(L^2(\R^d))$ and $F$ and $F_k$ are intertwining operators and thus bounded, the image of $\tilde{\nu}_k$ is contained in $\B(L^2(\R^d))$, as requested. \\
Next, it needs to be shown that $\tilde{\nu}_k$ is bounded: By the definition of $\| \cdot \|_{S_{i-1}}$, one has for every $\varphi \in CB(S_{i-1})$ 
\begin{eqnarray*}
\| \tilde{\nu}_k(\varphi)\|_{op}~=~\big\| F_k^* \circ F \circ \varphi \big(\pi_{\ell}^V \big) \circ F^* \circ F_k \big\|_{op} 
~\leq~ \big\|\varphi \big(\pi_{\ell}^V \big)\big\|_{op} 
~\leq~ \| \varphi\|_{S_{i-1}}.
\end{eqnarray*} 

In addition, one can easily observe that $\tilde{\nu}_k$ is involutive: For every $\varphi \in CB(S_{i-1})$
$$\tilde{\nu}_k(\varphi)^*~=~\big(F_k^* \circ F \circ \varphi \big(\pi_{\ell}^V \big) \circ F^* \circ F_k\big)^* 
~=~F_k^* \circ F \circ \varphi^* \big(\pi_{\ell}^V \big) \circ F^* \circ F_k ~=~\tilde{\nu}_k(\varphi^*).$$

Now, the last thing to check is the required convergence of Condition 3(b): For every $a \in C^*(G)$
\begin{eqnarray*}
\big\| \pi_k^V(a)- \tilde{\nu}_k\big(\F(a)_{|S_{i-1}}\big)\big\|_{op}&=&\big\| \pi_k^V(a)-F_k^* \circ F \circ \F(a)_{|S_{i-1}}\big(\pi_{\ell}^V \big) \circ F^* \circ F_k \big\|_{op} \\
&=&\big\| \pi_k^V(a)- F_k^* \circ F \circ \pi_{\ell}^V(a) \circ F^* \circ F_k\big\|_{op} \\
&=&\big\| F_k^* \circ \pi_{\ell_k}^V(a) \circ F_k- F_k^* \circ \pi(a) \circ F_k\big\|_{op} \\
&=& \big\| F_k^* \circ \big(\pi_{\ell_k}^V-\pi\big)(a) \circ F_k \big\|_{op} \\
&\leq& \big\|\pi_{\ell_k}^V(a)-\pi(a)\big\|_{op} \overset{k \to \infty}{\longrightarrow}0.
\end{eqnarray*}
Therefore, the representations $\big(\pi_{k}^V\big)_k$ and the constructed $(\tilde{\nu}_k)_k$ fulfill Condition 3(b) and thus, in this case, the claim is shown. 

\subsection{Second case}\label{second}
In the second case the situation that $\lambda=0$ or $c_j=0$ for every $j \in \{1,...,d\}$ must be considered. \\
In this case, 
\begin{eqnarray*}
 \langle{\ell_k},{[X^k_j,Y_j^k] }\rangle ~=~c^k_j \la_k \overset{k \to \infty}{\longrightarrow} c_j \lambda~=~0~~~\forall ~j \in \{1,...,d\},
 \end{eqnarray*}
while $c^k_j\la_k\ne 0$ for every $k$ and every $j \in \{1,...,d\}$. \\
Then $\ell_{| [\g,\g]}=0$ and so every limit orbit $O$ in the set $L((O_k)_k)$ has the dimension $0$. \\
~\\
As in Calculation (\ref{Pik-Formel1}) in Chapter \ref{Kapitel Pik-Formel}, identify $G$ again with $\R^{2d+2+p}$. From now on, this identification will be used most of the time. Only in some cases where one applies $\ell_k$ or $\ell$ and thus it is important to know whether one is using the basis depending on $k$ or the limit basis, the calculation will be done in the above defined bases $( \cdot)_k$ or $(\cdot)_{\infty}$.\\
~\\
Now, adapt the methods developed in \cite{ludwig-turowska} to this given situation. \\
Let $s=(s_1,...,s_d)$, $\al=(\al_1,\cdots, \al_d)$, $\be=(\be_1,\cdots, \be_d) \in \R^d$ and define
\begin{eqnarray*}
\eta_{k,\al,\be}(s)~=~
\eta_{k,\al,\be}(s_1,\cdots,s_d)
:=~e^{2 \pi i \al s}
\prod
\limits_{j=1}^d ~\big| \lambda_k c_j^k \big|^{\frac{1}{4}}~ \eta_j \bigg(| \lambda_k
c_j^k|^{\frac{1}{2}}\bigg(s_j + \frac{\be_j}{\lambda_kc_j^k}\bigg)\bigg).
\end{eqnarray*}

Moreover, let $c^k_{\al,\be} $ be the coefficient function defined by
\begin{eqnarray*}
 c^k_{\al,\be}(g):=\langle{\pi_k(g)\et_{k,\al,\be}},{\et_{k,\al,\be}}\rangle~~~\forall~g\in G \cong \R^{2d+2+p}
 \end{eqnarray*}
and $\ell_{\al,\be} $ the linear functional
$$\ell_{\al,\be}(g)~=~\ell_{\al,\be}(x,y,t,z,a)~:=~\al x+ \be y~~~\forall~g=(x,y,t,z,a) \in G \cong \R^{2d+2+p}.$$
Then, as in \cite{ludwig-turowska}, one can show by similar computations that the functions $c^k_{\al,\be} $ converge uniformly on compacta to the character $\ch_{\ell+\ell_{\al,\be}}$.

\subsubsection{Definition of the $\nu_k$'s}
For $0 \not= \tilde{\eta} \in L^2(G/P_k, \chi_{\ell_k}) \cong L^2(\R^d)$ let $$P_{\tilde{\eta}}:L^2(\R^d) \to \C \tilde{\eta},~ \xi \mapsto \tilde{\eta} \langle \xi, \tilde{\eta} \rangle.$$
Then $P_{\tilde{\eta}}$ is the orthogonal projection onto the space $\C \tilde{\eta}$. \\
Let $h \in C^*(G/U, \chi_{\ell})$. Again, identify $G/U$ with $\R^d \times \R^d \cong \R^{2d}$ and as already introduced in Chapter \ref{Kapitel Pik-Formel} in order to show the dependence on $k$, here the utilization of the limit basis will be expressed by an index $\infty$ if necessary:
$$h_{\infty}(x,y):=h((x,y)_{\infty}).$$ Now, $\hat{h}_{\infty}$ can be seen as a function in  $C_0(\ell + \u^{\perp}) \cong C_0(\R^{2d})$ and, using this identification, define the linear operator
$$\nu_k(h)~:=~\int \limits_{\R^{2d}}\hat{h}_{\infty}(\tilde{x}, \tilde{y}) ~P_{\eta_{k,\tilde{x}, \tilde{y}}}
~\frac{d(\tilde{x}, \tilde{y})}{\prod \limits_{j=1}^d \big| \lambda_k c_j^k \big|}.$$
~\\
Then, the following proposition holds: 
\begin{proposition}\label{nuk-prop}
~\\
\begin{enumerate}
\item For every $k \in \N$ and $h \in {\cal S}(G/U, \chi_{\ell})$ the integral defining $\nu_k(h)$ converges in the operator norm. 
\item The operator $\nu_k(h)$ is compact and $\| \nu_k(h) \|_{op} \klgl \| h \|_{C^*(G/U, \chi_{\ell})}$.
\item $\nu_k$ is involutive, i.e. $\nu_k(h)^*=\nu_k(h^*)$ for every $h \in C^*(G/U, \chi_{\ell})$. 
\end{enumerate}
\end{proposition}

Proof: 
\begin{enumerate}
\item Let $h \in \S(G/U, \chi_{\ell}) \cong \S(\R^{2d})$. Since 
\begin{eqnarray*}
\| P_{\eta_{k,\al,\be}} \|_{op}
~=~ \| \eta_{k,\al,\be} \|_2^2 
~=~1,
\end{eqnarray*}
one can estimate the operator norm of $\nu_k(h)$ as follows: 
\begin{eqnarray*}
\| \nu_k(h) \|_{op}
&=& \Bigg\|~ \int \limits_{\R^{2d}}\hat{h}_{\infty}(\tilde{x}, \tilde{y})~ P_{\eta_{k,\tilde{x}, \tilde{y}}}
~\frac{d(\tilde{x}, \tilde{y})}{\prod \limits_{j=1}^d \big| \lambda_k c_j^k \big|}~ \Bigg\|_{op}  \\
&\klgl& \int \limits_{\R^{2d}}\big|\hat{h}_{\infty}(\tilde{x}, \tilde{y})\big| ~\frac{d(\tilde{x}, \tilde{y})}{\prod \limits_{j=1}^d \big| \lambda_k c_j^k \big|}
~=~\frac{\big\| \hat{h} \big\|_{L^1(\R^{2d})}}{\prod \limits_{j=1}^d \big| \lambda_k c_j^k \big|}.
\end{eqnarray*}
Therefore, the convergence of the integral $\nu_k(h)$ in the operator norm is shown for \mbox{$h \in \S(\R^{2d}) \cong {\cal S}(G/U, \chi_{\ell})$}. 

\item First, let $h \in {\cal S}(G/U, \chi_{\ell}) \cong \S(\R^{2d})$. \\ 
Define for $s=(s_1,...,s_d) \in \R^d$
$$\eta_{k,\be}(s)~:=~\prod \limits_{j=1}^d ~\big| \lambda_k c_j^k \big|^{\frac{1}{4}}~ \eta_j \bigg(\big| \lambda_k
c_j^k \big|^{\frac{1}{2}}\bigg(s_j + \frac{\be_j}{\lambda_kc_j^k}\bigg)\bigg).$$
Then $$\eta_{k,\al,\be}(s)~=~e^{2 \pi i \al s} \eta_{k,\be}(s)$$ and thus one has for $\xi \in {\cal S}(\R^{d})$ and $s \in \R^d$ 
\begin{eqnarray}
&&\nn\nu_k(h) \xi(s) \\
\nn&&~~\\
\nn&=&\int \limits_{\R^{2d}}\hat{h}_{\infty}(\tilde{x}, \tilde{y}) ~\big\langle \xi, \eta_{k,\tilde{x}, \tilde{y}} \big\rangle \eta_{k,\tilde{x}, \tilde{y}}(s)~\frac{d(\tilde{x}, \tilde{y})}{\prod \limits_{j=1}^d \big| \lambda_k c_j^k \big|} \\
\nn&=& \int \limits_{\R^{2d}}\hat{h}_{\infty}(\tilde{x}, \tilde{y}) ~ \bigg(~\int \limits_{\R^d} \xi(r) \quer{\eta}_{k,\tilde{x}, \tilde{y}}(r)~dr \bigg) \eta_{k,\tilde{x}, \tilde{y}}(s)~\frac{d(\tilde{x}, \tilde{y})}{\prod \limits_{j=1}^d \big| \lambda_k c_j^k \big|}\\
\nn&=& \int \limits_{\R^{2d}}\hat{h}_{\infty}(\tilde{x}, \tilde{y}) ~ \bigg(~\int \limits_{\R^d} \xi(r)e^{-2 \pi i \tilde{x} r} \quer{\eta}_{k,\tilde{y}}(r)~dr \bigg) e^{2 \pi i \tilde{x} s}\eta_{k,\tilde{y}}(s)~\frac{d(\tilde{x}, \tilde{y})}{\prod \limits_{j=1}^d \big| \lambda_k c_j^k \big|} 
\end{eqnarray}
\begin{eqnarray}\label{Nuk-Formel}
\nn~~&=& \int \limits_{\R^{d}} \int \limits_{\R^{d}} \int \limits_{\R^{d}} \hat{h}_{\infty}(\tilde{x}, \tilde{y}) 
e^{2 \pi i \tilde{x}(s-r)}\xi(r) \quer{\eta}_{k,\tilde{y}}(r)~ dr ~\eta_{k,\tilde{y}}(s)~\frac{d\tilde{x}d \tilde{y}}{\prod \limits_{j=1}^d \big| \lambda_k c_j^k \big|}~~~~~~~~~~ \\
&=& \int \limits_{\R^{d}} \int \limits_{\R^{d}} \hat{h}_{\infty}^{2}(s-r, \tilde{y}) \xi(r) \quer{\eta}_{k,\tilde{y}}(r) \eta_{k,\tilde{y}}(s)~\frac{d\tilde{y}}{\prod \limits_{j=1}^d \big| \lambda_k c_j^k \big|}~dr.
\end{eqnarray}
Hence, as the kernel function 
$$h_K(s,r)~:=~\int \limits_{\R^{d}} \hat{h}_{\infty}^{2}(s-r, \tilde{y}) \quer{\eta}_{k,\tilde{y}}(r) \eta_{k,\tilde{y}}(s)~\frac{d\tilde{y}}{\prod \limits_{j=1}^d \big| \lambda_k c_j^k \big|}$$
of $\nu_k(h)$ is in ${\cal S}(\R^{2d})$, $\nu_k(h)$ is a compact operator. \\
Now it will be shown that $$\| \nu_k(h) \|_{op}~\klgl~\big\| \hat{h} \big\|_{\infty}:$$
For $\xi \in {\cal S}(\R^d)$ one has similar as in \cite{ludwig-turowska}
\begin{eqnarray}\label{estimnuh}
\nn&&\| \nu_k(h) \xi \|^2_2 \\
\nn&&\\
\nn&=&\int \limits_{\R^d} ~\Bigg|~ \int \limits_{\R^{d}} \int \limits_{\R^{d}} \hat{h}_{\infty}^{2}(s-r, \tilde{y}) \xi(r) \quer{\eta}_{k,\tilde{y}}(r) \eta_{k,\tilde{y}}(s)~\frac{d\tilde{y}}{\prod \limits_{j=1}^d \big| \lambda_k c_j^k \big|}~dr \Bigg|^2 ds \\
\nn&=&
\int \limits_{\R^d} ~\Bigg|~ \int \limits_{\R^{d}} \hat{h}_{\infty}^{2}(\cdot, \tilde{y}) \ast (\xi\quer{\eta}_{k,\tilde{y}})(s) \eta_{k,\tilde{y}}(s)~\frac{d\tilde{y}}{\prod \limits_{j=1}^d \big| \lambda_k c_j^k \big|}~\Bigg|^2 ds \\
\nn&\overset{\substack{Cauchy-\\Schwarz,}}{\underset{\| \eta_j \|_2=1}{\klgl}}&
\prod \limits_{j=1}^d \big| \lambda_k c_j^k \big| \int \limits_{\R^d}  \int \limits_{\R^{d}}\big| \hat{h}_{\infty}^{2}(\cdot, \tilde{y}) \ast (\xi\quer{\eta}_{k,\tilde{y}})(s)\big|^2 ds~\frac{d\tilde{y}}{\prod \limits_{j=1}^d \big| \lambda_k c_j^k \big|^2}~  \\
\nn&\overset{Plancherel}{\klgl}&
\big\| \hat{h} \big\|_{\infty}^2 \prod \limits_{j=1}^d \frac{1}{\big| \lambda_k c_j^k \big|}~ \int \limits_{\R^{d}} \| \xi\quer{\eta}_{k,\tilde{y}}\|_2^2~d\tilde{y} \\
\nn&\overset{\| \eta_j \|_2=1}{=}& \big\| \hat{h} \big\|_{\infty}^2  \| \xi \|_2^2. 
\end{eqnarray}
Thus, since ${\cal S}(\R^d)$ is dense in $L^2(\R^d)$,
$$\| \nu_k(h) \|_{op}~=~ \sup \limits_{\substack{\xi \in L^2(\R^d),\\ \| \xi \|_2=1}} \| \nu_k(h)(\xi) \|_2
~\klgl~\big\| \hat{h} \big\|_{\infty}$$
for $h \in \S(\R^{2d}) \cong {\cal S}(G/U, \chi_{\ell})$.
Therefore, with the density of ${\cal S}(G/U, \chi_{\ell})$ in $C^*(G/U, \chi_{\ell})$, one gets the compactness of the operator $\nu_k(h)$ for $h \in C^*(G/U, \chi_{\ell})$, as well as the desired inequality
$$\| \nu_k(h) \|_{op}~\klgl~\big\| \hat{h} \big\|_{\infty}~=~\| h \|_{C^*(G/U, \chi_{\ell})}.$$
\item The proof of the involutivity of $\nu_k$ is straightforward.
\end{enumerate}
\qed

~\\
This proposition firstly shows that the image of the operator $\nu_k$ is located in $\B(L^2(\R^d))=\B(\H_i)$ as required in Condition 3(b). Secondly, the proposition gives the boundedness and the involutivity of the linear mappings $\nu_k$ for every $k \in \N$. For the analysis of the sequence $(\pi_k)_k$, it remains to show the convergence condition. 

\subsubsection{Theorem - Second Case}
\begin{theorem}\label{theorem 2.fall}
~\\
Define as in Subsection \ref{CstUdef} 
\begin{eqnarray*}
p_{G/U}&:&L^1(G) \to L^1(G/U, \chi_{\ell}),~\\
&&p_{G/U}(f)(\tilde{g})~:=~\int \limits_U f(\tilde{g}u) \chi_{\ell}(u)~du~~~\forall ~\tilde{g} \in G ~~~\forall~ f \in L^1(G)
\end{eqnarray*}
and canonically extend $p_{G/U}$ to a mapping going from $C^*(G)$ to $C^*(G/U, \chi_{\ell})$. \\
Furthermore let $a \in C^*(G)$. Then 
$$\lim \limits_{k \to \infty} \big\| \pi_k(a) - \nu_k \big(p_{G/U}(a) \big) \big\|_{op}~=~0.$$
\end{theorem}

Proof: \\
For $u=(t,z,\dot{a},\ddot{a})_{\infty} \in U=span\{T,Z,A_1,...,A_{\tilde{p}},A_{\tilde{p}+1},...,A_p\}$ 
$$\chi_{\ell}(u)~=~e^{-2 \pi i \langle \ell, (t,z,\dot{a},\ddot{a})_{\infty} \rangle}~=~e^{-2 \pi i (t \rho+z \la)}$$
and therefore, identifying $U$ again with $\R \times \R \times \R^{\tilde{p}} \times \R^{p-\tilde{p}}$ and $L^1(G/U, \chi_{\ell})$ with $L^1(\R^{2d})$, for 
$f \in L^1(G) \cong L^1(\R^{2d+2+p})$ and $\tilde{g}=(\tilde{x},\tilde{y},0,0,0,0) \in \R^{2d}$ one has
\begin{eqnarray}\label{h gleich f-hut fall2}
\nn \big(p_{G/U}(f)\big)_{\infty}(\tilde{g})&=&\int \limits_{\R \times \R \times \R^{\tilde{p}} \times \R^{p-\tilde{p}}} f_{\infty}(\tilde{x},\tilde{y},\tilde{t},\tilde{z},\tilde{\dot{a}},\tilde{\ddot{a}})e^{-2 \pi i (\tilde{t} \rho+\tilde{z} \la)}d(0,0,\tilde{t}, \tilde{z}, \tilde{\dot{a}},\tilde{\ddot{a}}) \\
\nn&&~\\
&=&\hat{f}_{\infty}^{3,4,5,6}(\tilde{x},\tilde{y}, \rho, \la, 0,0),
\end{eqnarray}
whereat $f_{\infty}(\tilde{x},\tilde{y}, \rho, \la, 0,0)=f((\tilde{x},\tilde{y}, \rho, \la, 0,0)_{\infty})$. \\
~\\
Now, let $f \in {\cal S}(G) \cong \S(\R^{2d+2+p})$ such that its Fourier transform in $[\g,\g]$ has a compact support on $G \cong \R^{2d+2+p}$. If one then writes the elements $g$ of $G$ as $g=(x,y,t,z,\dot{a},\ddot{a})$ like above, whereat 
$$x \in span\{X_1\} \times...\times span \{X_d\},~~~y \in span\{Y_1\} \times...\times span\{Y_d\},~~~t \in span\{T\},$$ 
$$z \in span\{Z\},~~~\dot{a} \in span\{A_1\} \times ... \times span\{ A_{\tilde{p}}\},~~~\ddot{a} \in span\{A_{\tilde{p}+1}\} \times ... \times span\{ A_p\}$$
or, respectively 
$$x \in span\{X_1^k\} \times...\times span \{X_d^k\},~~~y \in span\{Y_1^k\} \times...\times span\{Y_d^k\},~~~t \in span\{T_k\},$$ 
$$z \in span\{Z_k\},~~~\dot{a} \in span\{A_1^k\} \times ... \times span\{ A^k_{\tilde{p}}\},~~~\ddot{a} \in span\{A_{\tilde{p}+1}^k\} \times ... \times span\{A_p^k\},$$
this means that the partial Fourier transform $\hat{f}^{4,5}$ has a compact support in $G$, since \mbox{$[\g,\g] =span\{Z_k, A_1^k,...,A_{\tilde{p}}^k\}=span\{ Z, A_1,...,A_{\tilde{p}}\}$}. \\
~\\
Moreover, let $s \in \R^d$ and define
$$\eta_{k,0}(s)~:=~\prod \limits_{j=1}^{d} \big|\lambda_k c_j^k \big|^{\frac{1}{4}} \eta_j \Big( ~\big|\lambda_k c_j^k\big|^{\frac{1}{2}}(s_j) \Big).$$
(Compare the definition of $\eta_{k, \beta}$ in the last proof.) \\
~\\
If $\xi\in\S(\R^d) $, one has
\begin{eqnarray*}
&&\big(\pi_k(f)- \nu_k \big(p_{G/U}(f) \big)\big)\xi (s) \\
&&~\\
&\overset{(\ref{Pik-Formel1}), (\ref{Nuk-Formel})}{=}&\int \limits_{\R^d} \hat{f}_k^{2,3,4,5,6}\Big(s-r,-\frac{\lambda_k c^k}{2}(s+r),\rho_k,\lambda_k,0,0\Big) \xi(r)~dr \\
&-& \int \limits_{\R^{d}} \int \limits_{\R^{d}} \widehat{p_{G/U}(f)}_{\infty}^{2}(s-r, \tilde{y}) \xi(r) \quer{\eta}_{k,\tilde{y}}(r) \eta_{k,\tilde{y}}(s)~\frac{d\tilde{y}}{\prod \limits_{j=1}^d \big| \lambda_k c_j^k \big|}~dr \\ 
&\overset{\| \eta_{k,0}\|_2=1}{\underset{(\ref{h gleich f-hut fall2})}{=}}&\int \limits_{\R^d} \int \limits_{\R^d} \hat{f}_k^{2,3,4,5,6}\Big(s-r,-\frac{\lambda_k c^k}{2}(s+r),\rho_k,\lambda_k,0,0\Big) \xi(r) \quer{\eta}_{k,0}(\tilde{y})\eta_{k,0}(\tilde{y})~d\tilde{y} dr \\
&-& \int \limits_{\R^{d}} \int \limits_{\R^{d}} \hat{f}_{\infty}^{2,3,4,5,6}(s-r, \tilde{y},\rho, \lambda,0,0) \xi(r) \quer{\eta}_{k,\tilde{y}}(r) \eta_{k,\tilde{y}}(s)~\frac{d\tilde{y}} {\prod \limits_{j=1}^d \big| \lambda_k c_j^k \big|}~dr \\
&=&\int \limits_{\R^d} \int \limits_{\R^d} \hat{f}_k^{2,3,4,5,6}\Big(s-r,-\frac{\lambda_k c^k}{2}(s+r),\rho_k,\lambda_k,0,0\Big) \xi(r) \quer{\eta}_{k,0}(\tilde{y})\eta_{k,0}(\tilde{y})~d \tilde{y} dr \\
&-& \int \limits_{\R^{d}} \int \limits_{\R^{d}} \hat{f}_{\infty}^{2,3,4,5,6}\big(s-r,\lambda_k c^k( \tilde{y}-s),\rho, \lambda,0,0\big) \xi(r)
\quer{\eta}_{k,0}(\tilde{y}+r-s) \eta_{k,0}(\tilde{y})~d\tilde{y} dr. \\
\end{eqnarray*}
The just obtained integrals are now divided into five parts. To do so, new functions $q_k$, $u_k$, $v_k$, $n_k$ and $w_k$ are defined:
\begin{eqnarray*}
q_k(s, \tilde{y})&:=&\int \limits_{\R^{d}} \xi(r)\quer{\eta}_{k,0}(\tilde{y}+r-s) 
\bigg(\hat{f}_k^{2,3,4,5,6}\Big(s-r,-\frac{\lambda_k c^k}{2}(s+r),\textcolor{red}{\rho_k},\lambda_k,0,0\Big) \\ 
&&~~~~~~~~~~~~~~~~~~~~~~~~~~~~~~~-~\hat{f}_k^{2,3,4,5,6}\Big(s-r,-\frac{\lambda_k c^k}{2}(s+r), \textcolor{red}{\rho},\lambda_k,0,0\Big)\bigg)dr,
\end{eqnarray*}
\begin{eqnarray*}
u_k(s, \tilde{y})&:=&\int \limits_{\R^{d}} \xi(r)\quer{\eta}_{k,0}(\tilde{y}+r-s) 
\bigg(\hat{f}_k^{2,3,4,5,6}\Big(s-r,-\frac{\lambda_k c^k}{2}(s+r),\rho,\textcolor{red}{\lambda_k},0,0\Big) \\  
&&~~~~~~~~~~~~~~~~~~~~~~~~~~~~~~~-~\hat{f}_k^{2,3,4,5,6}\Big(s-r,-\frac{\lambda_k c^k}{2}(s+r), \rho,\textcolor{red}{\lambda},0,0\Big)\bigg)dr,
\end{eqnarray*}
\begin{eqnarray*}
v_k(s, \tilde{y})&:=&\int \limits_{\R^{d}} \xi(r)\quer{\eta}_{k,0}(\tilde{y}+r-s) 
\bigg(\hat{f}_k^{2,3,4,5,6}\Big(s-r,\textcolor{red}{-\frac{\lambda_k c^k}{2}(s+r)}, \rho,\lambda,0,0\Big) \\
&&~~~~~~~~~~~~~~~~~~~~~~~~~~~~~~~-~ \hat{f}_k^{2,3,4,5,6}\big(s-r,\textcolor{red}{\lambda_k c^k(\tilde{y}-s)}, \rho,\lambda,0,0\big) \bigg)dr,
\end{eqnarray*}
\begin{eqnarray*}
n_k(s, \tilde{y})&:=&\int \limits_{\R^{d}} \xi(r)\quer{\eta}_{k,0}(\tilde{y}+r-s) 
\bigg(\textcolor{red}{\hat{f}_k^{2,3,4,5,6}}\big(s-r,\lambda_k c^k(\tilde{y}-s), \rho,\lambda,0,0\big) \\
&&~~~~~~~~~~~~~~~~~~~~~~~~~~~~~~~-~ \textcolor{red}{\hat{f}_{\infty}^{2,3,4,5,6}}\big(s-r,\lambda_k c^k(\tilde{y}-s), \rho,\lambda,0,0\big) \bigg)dr
\end{eqnarray*}
and
\begin{eqnarray*}
w_k(s)&:=&\int \limits_{\R^{d}} \int \limits_{\R^d} \xi(r)\eta_{k,0}(\tilde{y})\big(\quer{\eta}_{k,0}(\tilde{y})-\quer{\eta}_{k,0}(\tilde{y}+r-s)\big) \\ 
&&~~~~~~~\hat{f}_k^{2,3,4,5,6}\Big(s-r,-\frac{\lambda_k c^k}{2}(s+r),\rho_k,\lambda_k,0,0\Big)~d r d\tilde{y}.
\end{eqnarray*}
Then, 
\begin{eqnarray*}
\big(\pi_k(f)- \nu_k \big(p_{G/U}(f) \big)\big)\xi(s)
&=&\int \limits_{\R^d} q_k(s, \tilde{y}) \eta_{k,0}(\tilde{y})~d \tilde{y}~+~\int \limits_{\R^d} u_k(s, \tilde{y}) \eta_{k,0}(\tilde{y})~d \tilde{y} \\
&+&\int \limits_{\R^d} v_k(s, \tilde{y}) \eta_{k,0}(\tilde{y})~d \tilde{y}~+\int \limits_{\R^d} n_k(s, \tilde{y}) \eta_{k,0}(\tilde{y})~d \tilde{y}~+~w_k(s).
\end{eqnarray*}
In order to show that $$\big\|\pi_k(f)- \nu_k \big(p_{G/U}(f) \big) \big\|_{op} \overset{k \to \infty}{\longrightarrow} 0,$$ it suffices to prove that there are $\kappa_k$, $\gamma_k$, $\delta_k$, $\om_k$ and $\epsilon_k$ which are going to $0$ for $k \to \infty$, such that $$\|q_k\|_2\klgl\kappa_k \| \xi \|_2, ~~~\|u_k\|_2\klgl \gamma_k \| \xi \|_2, ~~~\|v_k\|_2 \klgl \delta_k \| \xi \|_2, ~~~\|n_k\|_2 \klgl \om_k \| \xi \|_2$$
$$\text{and}~~~\|w_k\|_2 \klgl \epsilon_k \| \xi \|_2.$$
~\\
First, regard the last factor of the function $q_k$:
\begin{eqnarray*}
&&\hat{f}_k^{2,3,4,5,6}\Big(s-r,-\frac{\lambda_k c^k}{2}(s+r),\rho_k,\lambda_k,0,0\Big) -\hat{f}_k^{2,3,4,5,6}\Big(s-r,-\frac{\lambda_k c^k}{2}(s+r), \rho,\lambda_k,0,0\Big) \\
&=&(\rho_k-\rho) \int \limits_0^1 \partial_3 \hat{f}_k^{2,3,4,5,6}\Big(s-r,-\frac{\lambda_k c^k}{2}(s+r),\rho+t(\rho_k-\rho), \lambda_k,0,0\Big)~dt.
\end{eqnarray*} 
Thus, since $f$ is a Schwartz function, one can find a constant $C_1>0$ (depending on $f$), such that
\begin{eqnarray*}
&&\bigg|~\hat{f}_k^{2,3,4,5,6}\Big(s-r,-\frac{\lambda_k c^k}{2}(s+r),\rho_k,\lambda_k,0,0\Big) -\hat{f}_k^{2,3,4,5,6}\Big(s-r,-\frac{\lambda_k c^k}{2}(s+r), \rho,\lambda_k,0,0\Big) ~\bigg| \\
&\klgl&|\rho_k-\rho|~ \frac{C_1}{(1+ \| s-r \|)^{2d}}.
\end{eqnarray*}
Hence, one gets the following estimation for $q_k$:
\begin{eqnarray*}
\|q_k \|_2^2&=& \int \limits_{\R^{2d}} | q_k(s, \tilde{y})|^2~d(s,\tilde{y}) \\
&\klgl& \int \limits_{\R^{2d}} \bigg( ~\int \limits_{\R^d} \big| \xi(r)\quer{\eta}_{k,0}(\tilde{y}+r-s)\big| | \rho_k-\rho|~
\frac{C_1}{(1+ \| s-r \|)^{2d}}~dr \bigg)^2 d(s,\tilde{y}) \\
&\overset{Cauchy-}{\underset{Schwarz}{\klgl}}& C_1^2 | \rho_k-\rho|^2 \int \limits_{\R^{2d}} \Bigg(~\int \limits_{\R^d}~
\bigg| \frac{\xi(r)\quer{\eta}_{k,0}(\tilde{y}+r-s)}{(1+ \| s-r \|)^d} \bigg|^2~dr \Bigg) \\ 
&&~~~~~~~~~~~~~~~~~~~~\Bigg(~\int \limits_{\R^d} \bigg(\frac{1}{(1+ \| s-r \|)^{d}} \bigg)^2~dr \Bigg) d(s, \tilde{y}) \\
&=& C_1' | \rho_k-\rho|^2 \int \limits_{\R^{3d}}
 \frac{ |\xi(r)|^2}{(1+ \| s-r \|)^{2d}}~ \big|\eta_{k,0}(\tilde{y}+r-s)\big|^2~d(r,s, \tilde{y}) \\
&\overset{\| \eta_{k,0}\|_2=1}{\klgl}& C_1'' | \rho_k-\rho|^2 \| \xi \|_2^2, 
\end{eqnarray*}
whereat $C_1'>0$ and $C_1''>0$ are matching constants depending on $f$. Thus, for $\kappa_k:=\sqrt{C_1''}~ | \rho_k-\rho|$, $\kappa_k \overset{k \to \infty}{\longrightarrow}0$, since $\rho_k \overset{k \to \infty}{\longrightarrow} \rho$, and 
$$\|q_k \|_2 \klgl \kappa_k \| \xi \|_2.$$

As $\lambda_k \overset{k \to \infty}{\longrightarrow} \lambda$, the estimation for the function $u_k$ can be done analogously. \\
~\\
Now, regard $v_k$. Like for $q_k$ and $u_k$, one has 
\begin{eqnarray*}
&&\hat{f}_k^{2,3,4,5,6}\Big(s-r,-\frac{\lambda_k c^k}{2}(s+r), \rho, \lambda,0,0\Big) - \hat{f}_k^{2,3,4,5,6}\big(s-r,\lambda_k c^k(\tilde{y}-s),\rho, \lambda,0,0\big) \\
&=& \lambda_k c^k \Big( ~\frac{1}{2}(r-s)- (r-s+ \tilde{y}) \Big) \\
&\cdot& \int \limits_0^1 \partial_2 \hat{f}_k^{2,3,4,5,6}\bigg(s-r, \lambda_k c^k(\tilde{y}-s)+t \lambda_k c^k \Big(~ \frac{1}{2}(r-s)- (r-s+ \tilde{y}) \Big), \rho, \lambda,0,0\bigg)~dt, 
\end{eqnarray*}
whereat $\cdot$ is the scalar product, and hence there exists again an on $f$ depending constant $C_3$ such that
\begin{eqnarray*}
&&\bigg|~\hat{f}_k^{2,3,4,5,6}\Big(s-r,-\frac{\lambda_k c^k}{2}(s+r), \rho, \lambda,0,0\Big) - \hat{f}_k^{2,3,4,5,6}\big(s-r,\lambda_k c^k(\tilde{y}-s),\rho, \lambda,0,0\big) \bigg| \\
&\klgl& | \lambda_k|~ \Big(\big\|c^k(r-s)\big\|+\big\|c^k(r-s+ \tilde{y})\big\|\Big)~\frac{C_3}{(1+\|s-r\|)^{2d+1}}.
\end{eqnarray*}
Therefore, defining $\tilde{\eta}_j(t):= \|t\| \eta_j(t)$, one gets a similar estimation for $v_k$:
\begin{eqnarray*}
&&\| v_k \|_2^2 \\
&&~\\
&=& \int \limits_{\R^{2d}} | v_k(s, \tilde{y})|^2~d(s,\tilde{y}) \\
&\klgl& \int \limits_{\R^{2d}} \bigg( ~\int \limits_{\R^d} \big| \xi(r)\quer{\eta}_{k,0}(\tilde{y}+r-s)\big|| \lambda_k |~ \Big(\big\|c^k(r-s)\big\|+\big\|c^k(r-s+ \tilde{y})\big\|\Big) \\
&&~~~~~~~~~~~~\frac{C_3}{(1+\|s-r\|)^{2d+1}} ~dr \bigg)^2 d(s,\tilde{y}) \\
&&\\
&\overset{Cauchy-}{\underset{Schwarz}{\klgl}}& C_3'\int \limits_{\R^{3d}} \frac{| \xi(r)|^2}{(1+\|s-r\|)^{2d+2}}~ \big|\eta_{k,0}(\tilde{y}+r-s)\big|^2| \lambda_k |^2 \\
&&~~~~~~~~\Big(\big\|c^k(r-s)\big\|+\big\|c^k(r-s+ \tilde{y})\big\|\Big)^2 ~d(r,s,\tilde{y}) \\
&&\\
&\klgl& 2C_3'\int \limits_{\R^{3d}} \frac{| \xi(r)|^2}{(1+\|s-r\|)^{2d+2}}~ \big|\eta_{k,0}(\tilde{y}+r-s)\big|^2| \lambda_k |^2 
\big\|c^k(r-s+ \tilde{y})\big\|^2 ~ d(r,s,\tilde{y}) \\
&+& 2C_3'\int \limits_{\R^{3d}} \frac{| \xi(r)|^2}{(1+\|s-r\|)^{2d+2}}~ \big|\eta_{k,0}(\tilde{y}+r-s)\big|^2| \lambda_k |^2 \big\|c^k(r-s)\big\|^2 ~ d(r,s,\tilde{y}) \\
&\klgl& 2C_3'| \lambda_k |^2 \Bigg(~\int \limits_{\R^{2d}} \frac{| \xi(r)|^2}{(1+\|s-r\|)^{2d+2}} ~d(r,s) \Bigg) \bigg(~\int \limits_{\R^d} \big\|c^k \tilde{y}\big\|^2\big|\eta_{k,0}(\tilde{y})\big|^2 d\tilde{y} \bigg) \\
&+& 2C_3' \big\| \lambda_k c^k \big\|^2 \int \limits_{\R^{3d}} \frac{| \xi(r)|^2}{(1+\|s-r\|)^{2d}}~\big|\eta_{k,0}(\tilde{y}+r-s)\big|^2 ~d(r,s,\tilde{y})\\
&\klgl& 2C_3' \big\| \lambda_k c^k \big\| \Bigg(~\int \limits_{\R^{2d}} \frac{| \xi(r)|^2}{(1+\|s-r\|)^{2d+2}} ~d(r,s) \Bigg) \\ 
&&\bigg(~\prod \limits_{j=1}^{d} \big|\lambda_k c_j^k\big|^{\frac{1}{2}} \int \limits_{\R} \Big\| \big|\lambda_k c_j^k\big|^{\frac{1}{2}} \tilde{y}_j\Big\|^2~\Big|\eta_j\Big(\big|\lambda_k c_j^k\big|^{\frac{1}{2}}(\tilde{y}_j)\Big)\Big|^2 d\tilde{y}_j \bigg) \\
&+& 2C_3' \big\| \lambda_k c^k\big\|^2 \int \limits_{\R^{3d}} \frac{|\xi(r)|^2}{(1+\|s-r\|)^{2d}}~\big|\eta_{k,0}(\tilde{y}+r-s)\big|^2 ~ d(r,s,\tilde{y})\\ 
&\overset{\| \eta_{k,0} \|_2=1}{\klgl}& 2C_3' \big\| \lambda_k c^k\big\|~ \| \xi \|_2^2 \bigg(~\prod \limits_{j=1}^{d} \int \limits_{\R} |\tilde{\eta}_j(\tilde{y}_j)|^2 d\tilde{y}_j \bigg) 
~+~ 2C_3' \big\| \lambda_k c^k\big\|^2 \| \xi \|_2^2 \\
&=& C_3'' \bigg(~\prod \limits_{j=1}^{d} \| \tilde{\eta}_j \|_2^2+ \big\| \lambda_k c^k \big\| \bigg) \big\|\lambda_k c^k \big\| ~\| \xi \|_2^2
\end{eqnarray*}
with constants $C_3'>0$ and $C_3''>0$, again depending on $f$. \\
Now, since $\lambda_k c^k \overset{k \to \infty}{\longrightarrow}0$, 
$\delta_k:=\bigg(C_3'' \bigg(~\prod \limits_{j=1}^{d} \| \tilde{\eta}_j \|_2^2+ \big\| \lambda_k c^k \big\| \bigg)\big\|\lambda_k c^k\big\|\bigg)^{\frac{1}{2}}$ fulfills $\delta_k \overset{k \to \infty}{\longrightarrow}0$ and 
$$\| v_k \|_2 ~\klgl~ \delta_k \| \xi \|_2.$$
~\\
For the estimation of $n_k$, the fact that the Fourier transform in $[\g,\g]$, i.e. in the 4th and 5th variable, $\hat{f}^{4,5}=:\tilde{f}$ has a compact support will be needed. Therefore, let the support of $f$ be located in the compact set 
$$K_1 \times K_2 \times K_3 \times K_4 \times K_5 \times K_6 \subset \R^{d} \times \R^{d}\times \R \times \R \times \R^{\tilde{p}} \times \R^{p- \tilde{p}}$$
and let $K:=K_2 \times K_3 \times K_6 \subset \R^{d}\times \R \times \R^{p-\tilde{p}}$. \\
Furthermore, since a Fourier transform is independent of the choice of the basis, the value of $\tilde{f}$ on $[\g,\g] =span\{Z_k, A_1^k,...,A_{\tilde{p}}^k\}=span\{ Z, A_1,...,A_{\tilde{p}}\}$ expressed in the $k$-basis and its value expressed in the limit basis are the same:
$$f(\cdot, \cdot, \cdot, (z,\dot{a})_k, \cdot)=f(\cdot, \cdot, \cdot, (z,\dot{a})_{\infty}, \cdot).$$
So, in the course of this proof, the limit basis will be chosen for the representation of the 4th and 5th position of an element $g$. Then
\begin{eqnarray*}
&&\hat{f}^{2,3,4,5,6}_k\big(s-r,\lambda_k c^k(\tilde{y}-s),\rho,\lambda,0,0\big) - \hat{f}_{\infty}^{2,3,4,5,6}\big(s-r,\lambda_k c^k(\tilde{y}-s),\rho,\lambda,0,0\big) \\
&=&\int \limits_{\R^{d} \times \R \times \R^{p-\tilde{p}}} \bigg(\tilde{f}_k(s-r,y,t,\la,0,\ddot{a})-\tilde{f}_{\infty}(s-r,y,t,\la,0,\ddot{a})\bigg)
e^{-2 \pi i(\la_k c^k(\tilde{y}-s)y+\rho t)} d(y,t,\ddot{a})\\
&=&\int \limits_{\R^{d} \times \R \times \R^{p-\tilde{p}}} \bigg(\tilde{f}((s-r,y,t)_k,(\la,0)_{\infty},(\ddot{a})_k)-\tilde{f}((s-r,y,t,\la,0,\ddot{a})_{\infty})\bigg) \\
&&~~~~~~~~~~~~~~e^{-2 \pi i(\la_k c^k(\tilde{y}-s)y+\rho t)} d(y,t,\ddot{a})\\
&=&\int \limits_{K} \bigg(\tilde{f}((s-r,y,t)_k,(\la,0)_{\infty},(\ddot{a})_k)-\tilde{f}((s-r,y,t,\la,0,\ddot{a})_{\infty})\bigg) 
e^{-2 \pi i(\la_k c^k(\tilde{y}-s)y+\rho t)} d(y,t,\ddot{a}).
\end{eqnarray*}

Furthermore, 
\begin{eqnarray*}
&&\tilde{f}((s-r,y,t)_k,(\la,0)_{\infty},(\ddot{a})_k)-\tilde{f}((s-r,y,t,\la,0,\ddot{a})_{\infty}) \\
&=& \tilde{f}\bigg(\sum \limits_{i=1}^d (s_i-r_i)X_i^k+\sum \limits_{i=1}^d y_i Y_i^k+tT_k+\la Z+\sum \limits_{i=\tilde{p}+1}^p a_i A_i^k \bigg) \\
&-& \tilde{f}\bigg(\sum \limits_{i=1}^d (s_i-r_i)X_i+\sum \limits_{i=1}^d y_i Y_i+tT+\la Z+\sum \limits_{i=\tilde{p}+1}^p a_i A_i \bigg) \\
&=& \bigg(\sum \limits_{i=1}^d (s_i-r_i)(X_i^k-X_i)+\sum \limits_{i=1}^d y_i (Y_i^k-Y_i)+t(T_k-T)+\sum \limits_{i=\tilde{p}+1}^p a_i (A_i^k-A_i)\bigg) \\
&\cdot& \int \limits_0^1 \partial \tilde{f}\bigg(\sum \limits_{i=1}^d (s_i-r_i)X_i+\sum \limits_{i=1}^d y_i Y_i+ tT+\la Z +\sum \limits_{i=\tilde{p}+1}^p a_i A_i \\
&&~~~~~~~~+ \tilde{t}\bigg(\sum \limits_{i=1}^d (s_i-r_i)(X_i^k-X_i)+\sum \limits_{i=1}^d y_i (Y_i^k-Y_i)+t(T_k-T)+\sum \limits_{i=\tilde{p}+1}^p a_i (A_i^k-A_i)\bigg)\bigg)d\tilde{t}.
\end{eqnarray*}
Since
$$(X_1^k,...,X_d^k,Y_1^k,...,Y_d^k,T_k,Z_k,A_1^k,...,A_p^k) \overset{k \to \infty}{\longrightarrow}(X_1,...,X_d,Y_1,...,Y_d,T,Z,A_1,...,A_p),$$
there exist $\om_k^i \overset{k \to \infty}{\longrightarrow}0$ for $i \in \{1,...,4\}$ and an on $f$ depending constant $C_4>0$ such that 
\begin{eqnarray*}
&&\bigg|~ \tilde{f}((s-r,y,t)_k,(\la,0)_{\infty},(\ddot{a})_k)-\tilde{f}((s-r,y,t,\la,0,\ddot{a})_{\infty}) \bigg| \\
&\leq&\Big(\| s-r\| \om_k^1+\|y\| \om_k^2+|t| \om_k^3+ \|\ddot{a}\| \om_k^4 \Big) ~ \frac{C_4}{(1+\|s-r\|)^{2d+1}}. 
\end{eqnarray*}

Now, with the help of the two calculations above, $\|n_k\|_2^2$ can be estimated:
\begin{eqnarray*}
&&\|n_k\|_2^2 \\
&&~\\
&=&\int \limits_{\R^{2d}}|n_k(s, \tilde{y})|^2 d(s, \tilde{y}) \\
&=&\int \limits_{\R^{2d}}\bigg|~ \int \limits_{\R^{d}}\xi(r)\quer{\eta}_{k,0}(\tilde{y}+r-s)
\bigg(\hat{f}^{2,3,4,5,6}_k\big(s-r,\lambda_k c^k(\tilde{y}-s),\rho,\lambda,0,0\big) \\
&&~~~~~~~~~~- \hat{f}_{\infty}^{2,3,4,5,6}\big(s-r,\lambda_k c^k(\tilde{y}-s),\rho,\lambda,0,0\big)\bigg)dr \bigg|^2 d(s, \tilde{y}) \\
&\leq& \int \limits_{\R^{2d}}\bigg( ~\int \limits_{\R^{d}} \big|\xi(r)\quer{\eta}_{k,0}(\tilde{y}+r-s)\big| \int \limits_{K} 
\Big(\| s-r\| \om_k^1+\|y\| \om_k^2+|t| \om_k^3+ \|\ddot{a}\| \om_k^4 \Big) \\
&&~~~~~~~~~~\frac{C_4}{(1+\|s-r\|)^{2d+1}}~
\big| e^{-2 \pi i(\la_k c^k(\tilde{y}-s)y+\rho t)}\big| d(y,t,\ddot{a}) dr \bigg)^2 d(s, \tilde{y}) \\
&=& C_4^2 \int \limits_{\R^{2d}} \bigg(~\int \limits_{\R^{d}} \frac{1}{(1+\|s-r\|)^{2d+1}} ~|\xi(r)| \big|\quer{\eta}_{k,0}(\tilde{y}+r-s)\big| \Big(\| s-r\| \om_k^5+\om_k^6 \Big) dr \bigg)^2 d(s, \tilde{y}) \\
&\underset{Schwarz}{\overset{Cauchy-}{\leq}}& C_4^2 \int \limits_{\R^{2d}} \bigg(~\int \limits_{\R^{d}} \frac{1}{(1+\|s-r\|)^{2d+2}} ~|\xi(r)|^2 \big|\quer{\eta}_{k,0}(\tilde{y}+r-s)\big|^2 \Big(\| s-r\| \om_k^5+\om_k^6 \Big)^2 dr \bigg) \\
&& ~~~~~~~~~\bigg(~\int \limits_{\R^{d}} \frac{1}{(1+\|s-r\|)^{2d}}~dr \bigg) d(s, \tilde{y}) \\
&\leq& C'_4\om_k^7 \int \limits_{\R^{2d}} \int \limits_{\R^{d}} \frac{1}{(1+\|s-r\|)^{2d}} ~|\xi(r)|^2 \big|\quer{\eta}_{k,0}(\tilde{y}+r-s)\big|^2 dr d(s, \tilde{y}) \\
&\overset{\| \eta_{k,0} \|_2=1}{=}& C''_4\om_k^7 \| \xi \|_2^2
\end{eqnarray*}
with constants $C'_4>0$ and $C''_4>0$ depending on $f$ and $\om_k^i \overset{k \to \infty}{\longrightarrow}0$ for $i \in \{5,...,7\}$. Thus, $\om_k:= \sqrt{C''_4 \om_k^7}$ fulfills $\om_k \overset{k \to \infty}{\longrightarrow} 0$ and 
$$\| n_k \|_2~\leq~\om_k \|\xi\|_2.$$
~\\
Last, it still remains to examine $w_k$:
\begin{eqnarray*}
&&\quer{\eta}_{k,0}(\tilde{y})-\quer{\eta}_{k,0}(\tilde{y}+r-s) \\
&=&\sum \limits_{j=1}^{d}(r_j-s_j) \int \limits_0^1 \partial_j \quer{\eta}_{k,0}(\tilde{y}+t(r-s))dt \\
&=&\sum \limits_{j=1}^d (r_j-s_j) \int \limits_0^1 \Bigg(~ \prod_{\substack{i=1\\i \not=j}}^d \big|\lambda_k c_i^k\big|^{\frac{1}{4}} \quer{\eta}_i\Big(\big|\lambda_k c_i^k\big|^{\frac{1}{2}}(\tilde{y}_i+t(r_i-s_i))\Big)\\
&&~~~~~~~~~~~~~~~~~~~~~~~~~~~~\big|\lambda_k c_j^k\big|^{\frac{3}{4}}
\partial \quer{\eta}_j\Big(\big|\lambda_k c_j^k\big|^{\frac{1}{2}}(\tilde{y}_j+t(r_j-s_j))\Big) \Bigg)dt. 
\end{eqnarray*}
Thus, since $f$ and the functions $(\eta_j)_{j \in \{1,...,d\}}$ are Schwartz functions, one can find an on $(\eta_j)_{j \in \{1,...,d\}}$ depending constant $C_5$ such that
\begin{eqnarray*}
&&\bigg|\Big(\quer{\eta}_{k,0}(\tilde{y})-\quer{\eta}_{k,0}(\tilde{y}+r-s)\Big) \hat{f}_k^{2,3,4,5,6}\Big(s-r, - \frac{\lambda_k c^k}{2}(s+r), \rho_k, \lambda_k,0,0\Big) \bigg| \\
&\klgl&\|r-s\|\bigg(\sum \limits_{j=1}^{d} \prod \limits_{i=1}^d \big|\lambda_k c_i^k\big|^{\frac{1}{4}}\big|\lambda_k c_j^k\big|^{\frac{1}{2}}\bigg)~\frac{C_5}{(1+\|r-s\|)^{2d+1}}.
\end{eqnarray*}
Now, one has the following estimation for $\|w_k \|_2$, which is again similar to the above ones:
\begin{eqnarray*}
\| w_k \|_2^2&=& \int \limits_{\R^d} |w_k(s)|^2 ds \\
&\klgl& \int \limits_{\R^d} \bigg(~ \int \limits_{\R^d} \int \limits_{\R^d} |\xi (r) | \big|\eta_{k,0}(\tilde{y})\big|~\|r-s\|\bigg(\sum \limits_{j=1}^{d} \prod \limits_{i=1}^d \big|\lambda_k c_i^k\big|^{\frac{1}{4}}\big|\lambda_k c_j^k\big|^{\frac{1}{2}}\bigg) \\
&&~~~~~~~~~~~~~~\frac{C_5}{(1+\|r-s\|)^{2d+1}}~dr d\tilde{y} \bigg)^2 ds \\
&\overset{Cauchy-}{\underset{Schwarz}{\klgl}}& C_5' \bigg(\sum \limits_{j=1}^{d} \prod \limits_{i=1}^d \big|\lambda_k c_i^k\big|^{\frac{1}{2}}\big|\lambda_k c_j^k\big|\bigg) \int \limits_{\R^{3d}} \frac{|\xi(r)|^2}{(1+\|r-s\|)^{2d+2}}~\big|\eta_{k,0}(\tilde{y})\big|^2 \|r-s\|^2~d(r,\tilde{y},s) \\
&\overset{\| \eta_{k,0}\|_2=1}{\klgl}& C_5'' \bigg(\sum \limits_{j=1}^{d} \prod \limits_{i=1}^d \big|\lambda_k c_i^k\big|^{\frac{1}{2}}\big|\lambda_k c_j^k\big|\bigg) \| \xi \|_2^2,
\end{eqnarray*}
whereat the constants $C_5'>0$ and $C_5''>0$ depend on $(\eta_j)_{j \in \{1,...,d\}}$. Therefore, for 
\mbox{$\epsilon_k:= \bigg(C_5'' \bigg(\sum \limits_{j=1}^{d} \prod \limits_{i=1}^d \big|\lambda_k c_i^k \big|^{\frac{1}{2}} \big|\lambda_k c_j^k \big|\bigg)\bigg)^{\frac{1}{2}}$}, the desired properties $\epsilon_k \overset{k \to \infty}{\longrightarrow}0$ and 
$$\| w_k \|_2~\klgl~\epsilon_k \| \xi \|_2$$
are fulfilled. \\
~\\
Thus, for those $f \in \S(\R^{2d+2+p}) \cong {\cal S}(G)$ whose Fourier transform in $[\g,\g]$ has a compact support,
$$ \big\|\pi_k(f)- \nu_k \big(p_{G/U}(f) \big) \big\|_{op}~=~\sup_{\substack{\xi \in L^2(\R^d)\\ \| \xi\|_2=1}} \big\|\big(\pi_k(f)- \nu_k \big(p_{G/U}(f) \big)\big)(\xi)\big\|_2
\overset{k \to \infty}{\longrightarrow}0.$$
Because of the density in $L^1(G)$ and thus in $C^*(G)$ of the set of Schwartz functions $f \in {\cal S}(G)$ whose partial Fourier transform has a compact support, the claim is true for general $a \in C^*(G)$. \\
\qed

\subsubsection{Transition to $\big(\pi_k^V \big)_k$}
As for every $k \in \N$ the two representations $\pi_k$ and $\pi_k^V$ are equivalent, there exist unitary intertwining operators
$$F_k: \H_{\pi_k^V} \cong L^2(\R^d) \to \H_{\pi_k} \cong L^2(\R^d)~~\text{with}~~F_k \circ \pi_k^V(a)=\pi_k(a) \circ F_k~~~\forall~ a \in C^*(G).$$
Futhermore, since the limit set $L \big( \big(\pi_k^V \big)_k\big)$ of the sequence $\big(\pi_k^V)_k$ is contained in $S_{i-1}$, as discussed in Section \ref{intro}, identifying $\wh{G}$ with the set of coadjoint orbits $\g^*/G$, one can restrict an operator field $\va \in CB(S_{i-1})$ to $L((O_k)_k)=\ell+\u^\perp $ and obtains an element in  $CB(\ell+\u^\perp) $. Thus, as
$$\big\{\F(a)_{| L((O_k)_k)}|~a \in C^*(G)\big\}=C_0(L((O_k)_k)) = C_0(\ell + \u^{\perp}),$$
one can define the $*$-isomorphism 
$$\tau:C_0(\R^{2d}) \cong C_0(\ell+\u^{\perp}) \to C^*(G/U, \chi_{\ell}) \cong C^*(\R^{2d}),~\F(a)_{| L((O_k)_k)} \mapsto p_{G/U}(a).$$
Now, define $\tilde{\nu}_k$ as 
$$\tilde{\nu}_k(\va)=F_k^* \circ (\nu_k \circ \tau)\big(\va_{| L((O_k)_k)}\big) \circ F_k~~~\forall~\va \in CB(S_{i-1}).$$
Since the image of $\nu_k$ is in $\B(L^2(\R^d))$ and $F_k$ is an intertwining operator and thus bounded, the image of $\tilde{\nu}_k$ is contained in $\B(L^2(\R^d))$ as well. \\
Moreover, the operator $\tilde{\nu}_k$ is bounded: From the boundedness of $\nu_k$ (Propostition \ref{nuk-prop}) and using that $\tau$ is an isomorphism, one gets for every $\va \in CB(S_{i-1})$
\begin{eqnarray*}
\| \tilde{\nu}_k(\va) \|_{op}&=& \big\|F_k^* \circ (\nu_k \circ \tau)\big(\va_{| L((O_k)_k)}\big) \circ F_k \big\|_{op} \\
&\leq& \big\|(\nu_k\circ \tau) \big(\va_{| L((O_k)_k)}\big)\big\|_{op} \\
&\leq& \big\|\tau \big(\va_{| L((O_k)_k)}\big)\big\|_{C^*(\R^{2d})} \\
&\leq&\big\|\big(\va_{| L((O_k)_k)}\big)\big\|_{\infty}~\leq~\|\va\|_{S_{i-1}}.
\end{eqnarray*}
The involutivity of $\tilde{\nu}_k$ follows from the involutivity of $\tau$ and $\nu_k$ (Proposition \ref{nuk-prop}). \\
Finally, the demanded convergence of Condition 3(b) can also be shown: With the above stated equivalence of the representations $\pi_k$ and $\pi_k^V$, one gets
\begin{eqnarray*}
\big\| \pi_k^V(a)- \tilde{\nu}_k\big(\F(a)_{|S_{i-1}}\big) \big\|_{op}&=& \big\| F_k^* \circ \pi_k(a) \circ F_k-F_k^* \circ (\nu_k \circ \tau)\big(\F(a)_{| L((O_k)_k)}\big) \circ F_k  \big\|_{op} \\
&=& \big\| F_k^* \circ \pi_k(a) \circ F_k - F_k^* \circ \nu_k \big(p_{G/U}(a)\big) \circ F_k  \big\|_{op} \\
&=& \big\| F_k^* \circ \big(\pi_k(a)-\nu_k \big(p_{G/U}(a) \big) \big) \circ F_k\big \|_{op} \\
&\leq& \big\| \nu_k \big(p_{G/U}(a) \big)-\pi_k(a)\big\|_{op} \overset{k \to \infty}{\longrightarrow}0.
\end{eqnarray*}
Therefore, the representations $\big(\pi_k^V \big)_k$ fulfill Property 3(b) and the conditions of the theorem are thus proved. 

\subsection{Third case}\label{third}
In the third and last case $\lambda \not=0$ and there exists $1\leq m<d $ such that $c_j \not=0$ for every $j \in \{1,...,m\}$ and $c_j=0$ for every 
$j \in \{m+1,...,d\}$. \\
This means that 
$$\langle{\ell_k},{[X^k_j,Y_j^k] }\rangle ~=~c^k_j \la_k \overset{k \to \infty}{\longrightarrow} c_j \lambda~=~0~~~\Longleftrightarrow~~~ j \in \{m+1,...,d\}.$$

In this case $\p:=\text{span}\{X_{m+1},...,X_d, Y_1,...,Y_d, T,Z,A_1,...,A_p\}$ is a polarization for $\ell$.  \\
Moreover, for $\tilde{\p}_k:=\text{span}\{X_{m+1}^k,...,X_d^k, Y_1^k,...,Y_d^k, T_k,Z_k, A_1^k,...,A_p^k\}$, one has $\tilde{\p}_k \overset{k \to \infty}{\longrightarrow}\p$. \\
Let $P:=exp(\p)$ and $\tilde{P}_k:=exp(\tilde{\p}_k)$. 

\subsubsection{Convergence of $(\pi_k)_k$ in $\wh{G}$}
Let 
$$(x)_{\infty}=(\dot{x},\ddot{x})_{\infty}~~\text{with}~~(\dot{x})_{\infty}:=(x_1,...,x_m)_{\infty}~~\text{and}~~ (\ddot{x})_{\infty}:=(x_{m+1},...,x_d)_{\infty}$$
and analogously 
$$(y)_{\infty}=(\dot{y},\ddot{y})_{\infty}~~\text{with}~~(\dot{y})_{\infty}:=(y_1,...,y_m)_{\infty}~~\text{and}~~(\ddot{y})_{\infty}:=(y_{m+1},...,y_d)_{\infty}.$$
\newpage
Moreover, as in Chapter \ref{definitions} above, let 
$$(a)_{\infty}=(\dot{a},\ddot{a})_{\infty}~~\text{with}~~(\dot{a})_{\infty}=(a_1,...,a_{\tilde{p}})_{\infty}~~\text{and}~~(\ddot{a})_{\infty}=(a_{\tilde{p}+1},...,a_p)_{\infty}$$
and let $$(g)_{\infty}=(\dot{x},\ddot{x},\dot{y},\ddot{y},t,z,\dot{a},\ddot{a})_{\infty}=(x,y,t,z,a)_{\infty}=(x,h)_{\infty}.$$ 
~\\
Now, let $\ddot{\alpha}:=(\alpha_{m+1},...,\al_d) \in \R^{d-m}$ and $\ddot{\be}:=(\be_{m+1},...,\be_d) \in \R^{d-m}$, consider $\ddot{\al}$ and $\ddot{\be}$ as elements of $\R^d$ identifying them with $(0,...,0,\alpha_{m+1},...,\al_d)$ and $(0,...,0,\be_{m+1},...,\be_d)$, respectively and let $$\tilde{\pi}:=\tilde{\pi}_{\ddot\al,\ddot\be}:=\ind_P^G \chi_{\ell+\ell_{\ddot{\al},\ddot{\be}}}.$$

Then, for a function $\dot{\xi}$ in the representation space $\H_{\tilde{\pi}}=L^2 \big(G/P, \chi_{\ell+\ell_{\ddot{\al},\ddot{\be}}}\big)$ of $\tilde{\pi}$ and 
for \mbox{$s_1,...,s_m \in \R$}, $(\dot{s})_{\infty}=(s_1,...,s_m)_{\infty}\in span \{X_1 \} \times...\times span \{X_m\}$ and $\dot{c}=(c_1,...,c_m)$, letting $\rho:=\langle \ell,T\rangle$ one has similarly as in (\ref{replk}):
\begin{eqnarray}\label{Pi-l-Formel}
 \nn\tilde{\pi}\big((g)_{\infty}\big) \dot{\xi}\big((\dot{s})_{\infty}\big)
~=~e^{2\pi i (-t \rho-z\la+ \la \dot{c}(((\dot{s})_{\infty}-\frac{1}{2}(\dot{x})_{\infty})(\dot{y})_{\infty}))}e^{-2 \pi i(\ddot{\al}(\ddot{x})_{\infty}+\ddot{\be}(\ddot{y})_{\infty})} \dot{\xi}\big((\dot{s}-\dot{x})_{\infty} \big),
\end{eqnarray}
since $\ell(Y_j)=\ell(X_j)=0$ for all $j \in \{1,...,d\}$. \\

From now on again, most of the time, $G$ will be identified with $\R^{2d+2+p}$.   
~\\
Define for $\ddot{s}=(s_{m+1},...,s_d) \in span \{X_{m+1}\} \times...\times span\{X_d\} \cong \R^{d-m}$
\begin{eqnarray*}
\ddot{\eta}_{k,\ddot{\al},\ddot{\be}}(\ddot{s})&:=&~e^{2 \pi i \ddot{\al} \ddot{s}}
\prod
\limits_{j=m+1}^d ~ \big| \lambda_k c_j^k \big|^{\frac{1}{4}}~ \eta_j \bigg(\big| \lambda_k
c_j^k \big|^{\frac{1}{2}}\bigg(s_j + \frac{\be_j}{\lambda_kc_j^k}\bigg)\bigg)
\end{eqnarray*}
and furthermore for $\dot{\xi} \in \H_{\tilde{\pi}}=L^2 \big(G/P, \chi_{\ell+\ell_{\ddot{\al},\ddot{\be}}} \big)\cong L^2(\R^m)$ and $s=(\dot{s},\ddot{s})$ in \mbox{$\big(span\{X_1\} \times...\times span\{X_m\}\big) \times \big(span\{X_{m+1}\} \times...\times span\{X_d\}\big) \cong \R^{m} \times \R^{d-m}$}
$$\xi_k(s):= \dot{\xi}(\dot{s}) \ddot{\eta}_{k,\ddot{\al},\ddot{\be}}(\ddot{s}).$$

Then, as above in the second case, the coefficient functions $c_{\ddot{\al},\ddot{\be}}^k$ defined by 
$$c_{\ddot{\al},\ddot{\be}}^k(g):=\langle \pi_k(g) \xi_k, \xi_k \rangle ~~~\forall ~g \in G \cong \R^{2d+2+p}$$
converge uniformly on compacta to $c_{\ddot{\al},\ddot{\be}}$ which in turn is defined by
$$c_{\ddot{\al},\ddot{\be}}(g):=\big\langle \tilde{\pi}_{\ddot{\al},\ddot{\be}}(g) \dot{\xi}, \dot{\xi} \big\rangle~~~\forall ~g \in G \cong \R^{2d+2+p}.$$

\subsubsection{Definition of the $\nu_k$'s} \label{nuk-def}
For $0 \not= \eta' \in L^2(\tilde{P}_k/P_k, \chi_{\ell_k}) \cong L^2(\R^{d-m})$ let
$$P_{\eta'}:L^2(\R^{d-m}) \to \C \eta',~ \xi \mapsto \eta' \langle \xi, \eta' \rangle.$$
Then $P_{\eta'}$ is the orthogonal projection onto the space $\C \eta'$. \\
Define now for $k \in \N$ and $h \in C^*(G/U, \chi_{\ell})$ the linear operator
$$\nu_k(h)~:=~\int \limits_{\R^{2(d-m)}} \pi_{\ell+ (\tilde{x},\tilde{y})}(h) \otimes P_{\ddot{\eta}_{k,\tilde{x}, \tilde{y}}}
~\frac{d(\tilde{x}, \tilde{y})}{\prod \limits_{j=m+1}^d \big| \lambda_k c_j^k \big|},$$
whereat $\pi_{\ell+ (\tilde{x},\tilde{y})}$ is defined as $ind_P^G \chi_{\ell+(\tilde{x},\tilde{y})}$ for an element $\ell+(\tilde{x},\tilde{y})$ located in  \mbox{$\ell+\big(\big(span\{X_{m+1}\}\times...\times span\{X_d\}\big) \times \big(span\{Y_{m+1}\}\times...\times span\{Y_d\}\big)\big)^* \cong \ell + \R^{2(d-m)}$.} \\
Thus, for $L^2(\R^d) \ni \xi= \sum \limits_{i=1}^{\infty} \dot{\xi}_i \otimes \ddot{\xi}_i$ with $\dot{\xi}_i \in L^2(\R^m)$ and $\ddot{\xi}_i \in L^2(\R^{d-m})$ for all $i \in \N$, one has 
$$\nu_k(h)(\xi)~:=~\sum \limits_{i=1}^{\infty} ~\int \limits_{\R^{2(d-m)}} \pi_{\ell+ (\tilde{x},\tilde{y})}(h)(\dot{\xi}_i) \otimes P_{\ddot{\eta}_{k,\tilde{x}, \tilde{y}}}(\ddot{\xi}_i)
~\frac{d(\tilde{x}, \tilde{y})}{\prod \limits_{j=m+1}^d \big| \lambda_k c_j^k \big|}.$$

\begin{proposition}
~\\
\begin{enumerate}
\item For every $k \in \N$ and $h \in {\cal S}(G/U, \chi_{\ell})$ the integral defining  $\nu_k(h)$ converges in the operator norm. 
\item The operator $\nu_k(h)$ is compact and $\| \nu_k(h) \|_{op} \klgl  \| h \|_{C^*(G/U, \chi_{\ell})}$.
\item $\nu_k$ is involutive. \\
\end{enumerate}
\end{proposition}

Proof: \\
Let $\K=\K(\l2{\R^m}) $ be the $C^* $-algebra of the compact operators on the Hilbert space $\l2{\R^{m}} $ and $C_0(\R^{2(d-m)}, \K) $ the $C^* $-algebra of all continuous mappings from $\R^{2(d-m)} $ into $\K $ vanishing at infinity. \\
Define for $\va \in C_0(\R^{2(d-m)},\K) $ and $k\in\N $ the linear operator 
\begin{eqnarray*}
 \mu_k(\va):=\int \limits_{\R^{2(d-m)}}\va(\tilde{x},\tilde{y})\ot P_{\ddot\et_{k,\tilde{x},\tilde{y}}}~\frac{d(\tilde{x},\tilde{y})}{\prod \limits_{j=m+1}^d \big| \lambda_k c_j^k \big|}
 \end{eqnarray*}
on $\l2{\R^d}$. Then, as $\F(h) \in C_0(\R^{2(d-m)},\K)$ for $h \in C^*(G/U, \chi_{\ell})$, 
$$\nu_k(h)=\mu_k(\F(h)).$$
\begin{enumerate}
\item
Since $\F(h) \in {\cal S}(\R^{2(d-m)},\K)$ for $h \in {\cal S}(G/U, \chi_{\ell})$ and since
\begin{eqnarray*}
 \noop{\mu_k(\va)}\leq \int \limits_{\R^{2(d-m)}}\noop{{\va(\tilde{x},\tilde{y})}}~\frac{d(\tilde{x},\tilde{y})}{\prod \limits_{j=m+1}^d \big| \lambda_k c_j^k \big|}
 \end{eqnarray*}
for every $\va \in \S(\R^{2(d-m)},\K) $ and $k\in\N $, the first assertion follows immediately.
\item
As $p_{G/U}$ is surjective from the space $\S(G)$ onto the space $\S(G/U, \chi_{\ell})$, for every \mbox{$h \in \S(G/U, \chi_{\ell}) \cong \S(\R^{2d})$} there exists a function $f \in \S(G) \cong \S(\R^{2d+2+p})$ such that \mbox{$h=p_{G/U}(f)$} and, as shown in the second case, for $\tilde{g}=(\dot{x}, \ddot{x},\dot{y}, \ddot{y},0,0,0,0) \in G/U \cong \R^{2d}$ one has
\begin{eqnarray}\label{h gleich f-hut}
h_{\infty}(\tilde{g})~=~\hat{f}_{\infty}^{5,6,7,8}(\dot{x}, \ddot{x},\dot{y},\ddot{y}, \rho, \la, 0,0),
\end{eqnarray}
where again $h_{\infty}=h((\cdot)_{\infty})$ and $f_{\infty}=f((\cdot)_{\infty})$. \\

Now, let $s_1,...,s_m \in \R$ and $\dot{s}=(s_1,...,s_m)_{\infty}= \sum \limits_{j=1}^m s_j X_j$ and moreover, let \mbox{$(g)_{\infty}=(\dot{x},\ddot{x},\dot{y},\ddot{y},t,z,\dot{a},\ddot{a})_{\infty}=(x,y,t,z,a)_{\infty}=(x,h)_{\infty}$} with $\dot{x},\ddot{x},\dot{y},\ddot{y},\dot{a},\ddot{a}$ and $\dot{c}$ as above.
Then, one gets for $\dot{\xi}_i \in L^2(\R^m)$
\begin{eqnarray*}
\pi_{\ell+ (\tilde{x},\tilde{y})}\big((g)_{\infty}\big) \dot{\xi}_i \big((\dot{s})_{\infty}\big)
~=~e^{2\pi i (-t \rho-z\la+ \la \dot{c}(((\dot{s})_{\infty}-\frac{1}{2}(\dot{x})_{\infty})(\dot{y})_{\infty})-\langle \tilde{x},(\ddot{x})_{\infty} \rangle- \langle\tilde{y},(\ddot{y})_{\infty}\rangle)} \dot{\xi}_i \big((\dot{s}-\dot{x})_{\infty}\big).
\end{eqnarray*}

Using the Equality (\ref{h gleich f-hut}) and identifying $G$ with $\R^{2d+2+p}$, one gets for a function $h = p_{G/U}(f) \in \S(G/U, \chi_{\ell}) \cong \S(\R^{2d})$ with $f \in \S(G) \cong \S(\R^{2d+2+p})$ \\
\begin{eqnarray}\label{pi-gleichung}
\nn\pi_{\ell+ (\tilde{x},\tilde{y})}(h) \dot{\xi}_i(\dot{s})&=& \int \limits_{\R^{2d}} \big(p_{G/U}(f) \big)_{\infty}(\tilde{g}) \pi_{\ell+ (\tilde{x},\tilde{y})}(\tilde{g}) \dot{\xi}_i(\dot{s})~d \tilde{g} \\
\nn&=&\int \limits_{\R^{2d}}~ \int \limits_{\R^{2+\tilde{p}+(p-\tilde{p})}} f_{\infty}(\tilde{g}u)\chi_{\ell}(u)~du~ \pi_{\ell+ (\tilde{x},\tilde{y})}(\tilde{g}) \dot{\xi}_i(\dot{s})~d \tilde{g} \\
\nn&=&\int \limits_{\R^{d} \times \R^d \times \R \times \R \times \R^{\tilde{p}} \times \R^{p-\tilde{p}}} f_{\infty}(\dot{x},\ddot{x},\dot{y}, \ddot{y},t,z,\dot{a},\ddot{a}) e^{2 \pi i( \la \dot{c}((\dot{s}-\frac{1}{2}\dot{x})\dot{y})-\langle \tilde{x},\ddot{x} \rangle- \langle\tilde{y},\ddot{y}\rangle)} \\
\nn&& ~~~~~~~~~~~~~~~~~~~~~~~~~~~e^{-2 \pi i (t \rho+z \la)} \dot{\xi}_i(\dot{s}-\dot{x}) ~d(\dot{x},\ddot{x},\dot{y},\ddot{y},t,z,\dot{a},\ddot{a}) \\
\nn&=&\int \limits_{\R^{2m}} \hat{f}_{\infty}^{2,4,5,6,7,8}(\dot{x},\tilde{x},\dot{y},\tilde{y},\rho,\la,0,0)e^{2\pi i \la \dot{c}((\dot{s}-\frac{1}{2}\dot{x})\dot{y})}\dot{\xi}_i(\dot{s}-\dot{x}) ~d(\dot{x},\dot{y}) \\
\nn&=&\int \limits_{\R^{2m}}\hat{h}_{\infty}^{2,4}(\dot{x},\tilde{x},\dot{y},\tilde{y})e^{2\pi i \la \dot{c}((\dot{s}-\frac{1}{2}\dot{x})\dot{y})}\dot{\xi}_i(\dot{s}-\dot{x}) ~d(\dot{x},\dot{y}) \\
&=& \int \limits_{\R^m}\hat{h}_{\infty}^{2,3,4}\Big(\dot{x},\tilde{x},\la \dot{c}\Big(~\frac{1}{2}\dot{x}-\dot{s}\Big),\tilde{y}\Big) \dot{\xi}_i(\dot{s}-\dot{x}) ~d\dot{x}. \label{pi-gleichung 2}
\end{eqnarray}

Regard now the second factor $P_{\ddot{\eta}_{k,\tilde{x}, \tilde{y}}}$ of the tensor product: \\
As in the second case above, define $$\ddot{\eta}_{k,\ddot{\be}}(\ddot{s})~:=~\prod \limits_{j=m+1}^d ~\big| \lambda_k c_j^k \big|^{\frac{1}{4}}~ \eta_j \bigg(\big| \lambda_k c_j^k \big|^{\frac{1}{2}}\bigg(s_j + \frac{\be_j}{\lambda_kc_j^k}\bigg)\bigg).$$ 
Then $$\ddot{\eta}_{k,\ddot{\al}, \ddot{\be}}(\ddot{s})=e^{2 \pi i \ddot{\al}\ddot{s}} \ddot{\eta}_{k,\ddot{\be}}(\ddot{s})$$ and therefore with $\ddot{\xi}_i \in L^2(\R^{d-m})$
\begin{eqnarray*}
P_{\ddot{\eta}_{k,\tilde{x}, \tilde{y}}}(\ddot{\xi}_i)(\ddot{s})
&=& \big\langle \ddot{\xi}_i,\ddot{\eta}_{k,\tilde{x}, \tilde{y}} \big\rangle \ddot{\eta}_{k,\tilde{x}, \tilde{y}}(\ddot{s}) \\
&=& \bigg(~\int \limits_{\R^{d-m}} \ddot{\xi}_i(\ddot{r}) \quer{\ddot{\eta}}_{k,\tilde{x}, \tilde{y}}(\ddot{r})~d\ddot{r} \bigg) \ddot{\eta}_{k,\tilde{x}, \tilde{y}}(\ddot{s}) \\
&=& \bigg(~\int \limits_{\R^{d-m}} \ddot{\xi}_i(\ddot{r})e^{-2 \pi i \tilde{x}\ddot{r}} ~\quer{\ddot{\eta}}_{k,\tilde{y}}(\ddot{r})~d\ddot{r} \bigg) e^{2 \pi i \tilde{x}\ddot{s}}\ddot{\eta}_{k,\tilde{y}}(\ddot{s}) \\
&=& \int \limits_{\R^{d-m}} \ddot{\xi}_i(\ddot{r})e^{2 \pi i \tilde{x}(\ddot{s}-\ddot{r})} \quer{\ddot{\eta}}_{k,\tilde{y}}(\ddot{r})~d\ddot{r}~ \ddot{\eta}_{k,\tilde{y}}(\ddot{s}). \\
\end{eqnarray*} 
Joining together the calculation above and the one for the first factor of the tensor product (\ref{pi-gleichung 2}), one gets for $\S(\R^d) \ni \xi= \sum \limits_{i=1}^{\infty} \dot{\xi}_i \otimes \ddot{\xi}_i$ with $\dot{\xi}_i \in \S(\R^m)$ and $\ddot{\xi}_i \in \S(\R^{d-m})$ for all $i \in \N$, $h \in \S(\R^{2d})$ and $s=(\dot{s},\ddot{s}) \in \R^d$
\begin{eqnarray}\label{nuk-formel2}
\nn&&\nu_k(h)(\xi)(s) \\
\nn&&\\
\nn&=&\sum \limits_{i=1}^{\infty} ~\int \limits_{\R^{2(d-m)}} \pi_{\ell+ (\tilde{x},\tilde{y})}(h)(\dot{\xi}_i)(\dot{s}) \cdot P_{\ddot{\eta}_{k,\tilde{x}, \tilde{y}}}(\ddot{\xi}_i)(\ddot{s})~\frac{d(\tilde{x}, \tilde{y})}{\prod \limits_{j=m+1}^d \big| \lambda_k c_j^k \big|} \\
\nn&=& \sum \limits_{i=1}^{\infty} ~\int \limits_{\R^{2(d-m)}}  
\bigg(~\int \limits_{\R^m} \hat{h}_{\infty}^{2,3,4}\Big(\dot{x},\tilde{x},\la \dot{c} \Big(~\frac{1}{2}\dot{x}-\dot{s}\Big),\tilde{y}\Big) \dot{\xi}_i(\dot{s}-\dot{x}) ~d\dot{x}\bigg) \\
\nn& \cdot& \bigg(~\int \limits_{\R^{d-m}} \ddot{\xi}_i(\ddot{r})e^{2 \pi i \tilde{x}(\ddot{s}-\ddot{r})} \quer{\ddot{\eta}}_{k,\tilde{y}}(\ddot{r})~d\ddot{r}~ \ddot{\eta}_{k,\tilde{y}}(\ddot{s})\bigg)~\frac{d(\tilde{x}, \tilde{y})}{\prod \limits_{j=m+1}^d \big| \lambda_k c_j^k \big|} \\
\nn&=& \sum \limits_{i=1}^{\infty} ~\int \limits_{\R^{2(d-m)}} \int \limits_{\R^{d-m}} \int \limits_{\R^{m}} \hat{h}_{\infty}^{2,3,4}\Big(\dot{x},\tilde{x},\la \dot{c}\Big(~\frac{1}{2}\dot{x}-\dot{s}\Big),\tilde{y}\Big) \dot{\xi}_i(\dot{s}-\dot{x}) \\
\nn&&~~~~~~~~~~~~~~~~~~~~~~~~~~~\ddot{\xi}_i(\ddot{r})e^{2 \pi i \tilde{x}(\ddot{s}-\ddot{r})} \quer{\ddot{\eta}}_{k,\tilde{y}}(\ddot{r})~d\dot{x} d\ddot{r}~ \ddot{\eta}_{k,\tilde{y}}(\ddot{s})~\frac{d(\tilde{x}, \tilde{y})}{\prod \limits_{j=m+1}^d \big| \lambda_k c_j^k \big|} \\
\nn&=& \sum \limits_{i=1}^{\infty} ~\int \limits_{\R^{d-m}} \int \limits_{\R^{d-m}} \int \limits_{\R^{m}} \hat{h}_{\infty}^{3,4}\Big(\dot{x},\ddot{s}-\ddot{r},\la \dot{c}\Big(~\frac{1}{2}\dot{x}-\dot{s}\Big),\tilde{y}\Big) \dot{\xi}_i(\dot{s}-\dot{x}) \\
\nn&&~~~~~~~~~~~~~~~~~~~~~~~~\ddot{\xi}_i(\ddot{r}) \quer{\ddot{\eta}}_{k,\tilde{y}}(\ddot{r})~d\dot{x} d\ddot{r}~ \ddot{\eta}_{k,\tilde{y}}(\ddot{s})~\frac{d\tilde{y}}{\prod \limits_{j=m+1}^d \big| \lambda_k c_j^k \big|} \\
\nn&=& \sum \limits_{i=1}^{\infty} ~\int \limits_{\R^{d-m}} \int \limits_{\R^{d-m}} \int \limits_{\R^{m}} \hat{h}_{\infty}^{3,4}\Big(\dot{s}-\dot{x},\ddot{s}-\ddot{r},\frac{\la}{2} \dot{c}(-\dot{x}-\dot{s}),\tilde{y} \Big) \dot{\xi}_i(\dot{x}) \\
\nn&&~~~~~~~~~~~~~~~~~~~~~~~~\ddot{\xi}_i(\ddot{r}) \quer{\ddot{\eta}}_{k,\tilde{y}}(\ddot{r})~d\dot{x} d\ddot{r}~ \ddot{\eta}_{k,\tilde{y}}(\ddot{s})~\frac{d\tilde{y}}{\prod \limits_{j=m+1}^d \big| \lambda_k c_j^k \big|} \\
\nn&=& \int \limits_{\R^{d}} \int \limits_{\R^{d-m}} \hat{h}_{\infty}^{3,4}\Big(\dot{s}-\dot{x},\ddot{s}-\ddot{r},\frac{\la}{2} \dot{c}(-\dot{x}-\dot{s}),\tilde{y}\Big)  
\quer{\ddot{\eta}}_{k,\tilde{y}}(\ddot{r})~ \ddot{\eta}_{k,\tilde{y}}(\ddot{s})~\frac{d\tilde{y}}{\prod \limits_{j=m+1}^d \big| \lambda_k c_j^k \big|} \\
\label{nuk-formel2.2}&&~~~~~~~~~~\xi(\dot{x},\ddot{r})~d(\dot{x},\ddot{r}).
\end{eqnarray}

Therefore, the kernel function 
$$h_K((\dot s,\ddot s),(\dot{x},\ddot{r}))~:=~\int \limits_{\R^{d-m}}\hat{h}_{\infty}^{3,4}\Big(\dot{s}-\dot{x},\ddot{s}-\ddot{r},\frac{\la}{2} \dot{c}(-\dot{x}-\dot{s}),\tilde{y}\Big)  
\quer{\ddot{\eta}}_{k,\tilde{y}}(\ddot{r})~ \ddot{\eta}_{k,\tilde{y}}(\ddot{s})~\frac{d\tilde{y}}{\prod \limits_{j=m+1}^d \big| \lambda_k c_j^k \big|}$$
of $\nu_k(h)$ is contained in ${\cal S}(\R^{2d})$ and thus $\nu_k(h)$ is a compact operator for $h \in \S(\R^{2d}) \cong {\cal S}(G/U, \chi_{\ell})$ and with the density of ${\cal S}(G/U, \chi_{\ell})$ in $C^*(G/U, \chi_{\ell})$, it is compact for every $h \in C^*(G/U, \chi_{\ell})$.\\

Now, it is shown that for every $\va \in C_0(\R^{2(d-m)},\K)$
\begin{eqnarray*}
 \noop{\mu_k(\va)}~\leq ~\no{\va}_{\infty}~:=~\sup_{(\tilde{x},\tilde{y})\in\R^{2(d-m)}}\noop{\va(\tilde{x},\tilde{y})}.
 \end{eqnarray*}
For this, for any $\psi\in\l2{\R^d}$, define 
$$f_{\psi,k}(\tilde{x},\tilde{y})(\dot{s}):=\int \limits_{\R^{d-m}}\psi(\dot{s},\ddot{s})\ol{\ddot \et}_{k,\tilde{x},\tilde{y}}(\ddot{s})~d \ddot{s} ~~~\forall~(\tilde{x},\tilde{y}) \in \R^{2(d-m)}~~\forall~ \dot s\in \R^m.$$
Then, as 
\begin{eqnarray*}
 \Id_{\l2{\R^{d-m}}}~=~\int \limits_{\R^{2(d-m)}}P_{\ddot\et_{k,\tilde{x},\tilde{y}}}~\frac{d(\tilde{x},\tilde{y})}{\prod \limits_{j=m+1}^d \big| \lambda_k c_j^k \big|}, 
 \end{eqnarray*}
one gets the identity 
\begin{eqnarray}\label{puid}
\no\psi_{\l2{\R^d}}^2~=~ \int \limits_{\R^{2(d-m)}}\| f_{\psi,k}(\tilde{x},\tilde{y})\|_{L^2(\R^m)}^2~\frac{d(\tilde{x},\tilde{y})}{\prod \limits_{j=m+1}^d \big| \lambda_k c_j^k \big|}.
\end{eqnarray}
Now, for $\xi,\psi \in  \l2{\R^d}$ 
\begin{eqnarray*}
 & &\big| \langle{\mu_k(\va)\xi},{\psi}\rangle_{\l2{\R^d}} \big|\\
 &&~\\
 &=&\bigg\vert~\int \limits_{\R^{2(d-m)}}\langle{(\va(\tilde{x},\tilde{y})\ot P_{\ddot \et_{k,\tilde{x},\tilde{y}}})\xi},{\psi}\rangle_{\l2{\R^d}}~\frac{d(\tilde{x},\tilde{y})}{\prod \limits_{j=m+1}^d \big| \lambda_k c_j^k \big|} \bigg\vert\\
 &=& \bigg\vert~
 \int \limits_{\R^{2(d-m)}}\langle{(\va(\tilde{x},\tilde{y})\ot\Id_{\l2{\R^{d-m}}})\ci(\Id_{\l2{\R^m}}\ot P_{\ddot \et_{k,\tilde{x},\tilde{y}}})\xi},{(\Id_{\l2{\R^m}}\ot P_{\ddot \et_{k,\tilde{x},\tilde{y}}})\psi}\rangle_{\l2{\R^d}} \\
 &&~~~~~~~~~~~~~\frac{d(\tilde{x},\tilde{y})}{\prod \limits_{j=m+1}^d \big| \lambda_k c_j^k \big|} \bigg\vert\\
 &=&
   \bigg\vert~\int \limits_{\R^{2(d-m)}}\langle{\va(\tilde{x},\tilde{y})f_{\xi,k}(\tilde{x},\tilde{y})},{f_{\psi,k}(\tilde{x},\tilde{y})}\rangle_{\l2{\R^m}}~\frac{d(\tilde{x},\tilde{y})}{\prod \limits_{j=m+1}^d \big| \lambda_k c_j^k \big|} \bigg\vert\\
  &\overset{Cauchy-}{\underset{Schwarz}{\leq}}& 
 \bigg(~\int \limits_{\R^{2(d-m)}}\no{\va(\tilde{x},\tilde{y})f_{\xi,k}(\tilde{x},\tilde{y})}^2_{\l2{\R^m}}\frac{d(\tilde{x},\tilde{y})}{\prod \limits_{j=m+1}^d \big| \lambda_k c_j^k \big|}\bigg)^{\frac{1}{2}} \\
 && \bigg(~\int \limits_{\R^{2(d-m)}}\no{f_{\psi,k}(\tilde{x},\tilde{y})}^2_{\l2{\R^m}}\frac{d(\tilde{x},\tilde{y})}{\prod \limits_{j=m+1}^d \big| \lambda_k c_j^k \big|} \bigg)^{\frac{1}{2}} \\
  &{\overset{(\ref{puid})}{\leq}}&\sup_{(\tilde{x},\tilde{y})\in\R^{2(d-m)}}\noop{{\va(\tilde{x},\tilde{y})}}\bigg(~\int \limits_{\R^{2(d-m)}}\no{f_{\xi,k}(\tilde{x},\tilde{y})}^2_{\l2{\R^m}}\frac{d(\tilde{x},\tilde{y})}{\prod \limits_{j=m+1}^d \big| \lambda_k c_j^k \big|}\bigg)^{\frac{1}{2}} \no{\psi}_{\l2{\R^d}}\\
  &\overset{(\ref{puid})}{\leq}& \no\va_{\infty}\no{\xi}_{\l2{\R^d}}\no {\psi}_{\l2{\R^d}}.
 \end{eqnarray*}
Hence, for every $h \in C^*(G/U, \chi_{\ell})$,
$$\noop{\nu_k(h)}~=~\noop{\mu_k(\F(h))}~\leq ~\no{\F(h)}_{\infty}~=~\|h\|_{C^*(G/U, \chi_{\ell})}. $$
\item To show that $\nu_k$ is involutive is as straightforward as in the second case. \\
\end{enumerate}
\qed
\newpage

The demanded convergence of Condition 3(b) remains to be shown: 

\subsubsection{Theorem - Third Case}
\begin{theorem}
~\\
For $a \in C^*(G)$
$$\lim \limits_{k \to \infty} \big\| \pi_k(a) - \nu_k\big(p_{G/U}(a)\big) \big\|_{op}~=~0.$$
\end{theorem}

Proof: \\
Let $f \in {\cal S}(G) \cong \S(\R^{2d+2+p})$ such that its Fourier transform in $[\g,\g]$ has a compact support on $G \cong \R^{2d+2+p}$. In the setting of this third case, this means that $\hat{f}^{6,7}$ has a compact support in $G$ (see Theorem \ref{theorem 2.fall}). \\
~\\
Now, identify $G$ with $\R^{2d+2+p}$ again, let $\xi \in L^2(\R^d)$ and $s=(s_1,...,s_d)=(\dot{s},\ddot{s})$ be located in $\R^m \times \R^{d-m} \cong \R^d$ and define
$$\ddot{\eta}_{k,0}(\ddot{s})~:=~\prod \limits_{j=m+1}^{d} \big|\lambda_k c_j^k \big|^{\frac{1}{4}} \eta_j \Big(\big|\lambda_k c_j^k \big|^{\frac{1}{2}}(s_j)\Big).$$
Moreover, let $\dot{c}=(c_1,...,c_m)$, $\ddot{c}=(c_{m+1},...,c_d)=(0,...,0)$, $\dot{c}^k=(c_1^k,...,c_m^k)$ and \mbox{$\ddot{c}^k=(c_{m+1}^k,...,c_d^k)$}. \\
As in the second case, the expression $\big(\pi_k(f)- \nu_k\big(p_{G/U}(f)\big)\big)$ is now going to be regarded, composed into several parts and then estimated:
For this, Equation (\ref{Pik-Formel1}) from Chapter \ref{Kapitel Pik-Formel} will be used again but its notation needs to be adapted: 
\begin{eqnarray*}
\pi_k(f)(s)&\overset{(\ref{Pik-Formel1})}{=}&\int \limits_{\R^d} \hat f_k^{2,3,4,5,6}\Big(s-r,-\frac{\lambda_k c^k}{2}(s+r),\rho_k,\lambda_k,0,0 \Big) \xi(r)~dr \\
&=&\int \limits_{\R^d} \hat f_k^{3,4,5,6,7,8} \Big(\dot{s}-\dot{r},\ddot{s}-\ddot{r},-\frac{\lambda_k \dot{c}^k}{2}(\dot{s}+\dot{r}),-\frac{\lambda_k \ddot{c}^k}{2}(\ddot{s}+\ddot{r}),\rho_k,\lambda_k,0,0 \Big) \xi(\dot{r},\ddot{r})~d(\dot{r},\ddot{r}). \\
\end{eqnarray*}

Using the above equation, (\ref{nuk-formel2.2}) and the fact that $p_{G/U}(f)=\hat{f}^{5,6,7,8}(\cdot, \cdot, \cdot, \cdot, \rho, \la,0,0)$, one gets
\begin{eqnarray*}
&&\big(\pi_k(f)- \nu_k\big(p_{G/U}(f)\big)\big)\xi(s) \\
&&~\\
&\overset{(\ref{nuk-formel2.2})}{=}&\int \limits_{\R^d} \hat f_k^{3,4,5,6,7,8}\Big(\dot{s}-\dot{r},\ddot{s}-\ddot{r},-\frac{\lambda_k \dot{c}^k}{2}(\dot{s}+\dot{r}),-\frac{\lambda_k \ddot{c}^k}{2}(\ddot{s}+\ddot{r}),\rho_k,\lambda_k,0,0\Big) \xi(\dot{r},\ddot{r})~d(\dot{r},\ddot{r}) \\
&-& \int \limits_{\R^{d}} \int \limits_{\R^{d-m}} \widehat{p_{G/U}(f)}_{\infty}^{3,4}\Big(\dot{s}-\dot{r},\ddot{s}-\ddot{r},\frac{\la}{2} \dot{c}(-\dot{r}-\dot{s}),\tilde{y}\Big)  
\quer{\ddot{\eta}}_{k,\tilde{y}}(\ddot{r}) \ddot{\eta}_{k,\tilde{y}}(\ddot{s})~\frac{d\tilde{y}}{\prod \limits_{j=m+1}^d \big| \lambda_k c_j^k \big|} \\
&&~~~~\xi(\dot{r},\ddot{r})~d(\dot{r},\ddot{r}) \\
&&~~\\
&\overset{\| \tilde{\eta}_{k,0}\|_2=1}{=}& \int \limits_{\R^d} \int \limits_{\R^{d-m}} \hat f_k^{3,4,5,6,7,8}\Big(\dot{s}-\dot{r},\ddot{s}-\ddot{r},-\frac{\lambda_k \dot{c}^k}{2}(\dot{s}+\dot{r}),-\frac{\lambda_k \ddot{c}^k}{2}(\ddot{s}+\ddot{r}),\rho_k,\lambda_k,0,0\Big) \\
&&~~~~~~~~~~~\quer{\ddot{\eta}}_{k,0}(\tilde{y})\ddot{\eta}_{k,0}(\tilde{y})~d\tilde{y}~\xi(\dot{r},\ddot{r})~d(\dot{r},\ddot{r}) \\ 
&&~~\\
&-& \int \limits_{\R^{d}} \int \limits_{\R^{d-m}} \hat{f}_{\infty}^{3,4,5,6,7,8}\Big(\dot{s}-\dot{r},\ddot{s}-\ddot{r},\frac{\la}{2} \dot{c}(-\dot{r}-\dot{s}),\tilde{y},\rho,\lambda,0,0\Big)  
\quer{\ddot{\eta}}_{k,\tilde{y}}(\ddot{r}) \ddot{\eta}_{k,\tilde{y}}(\ddot{s})~\frac{d\tilde{y}}{\prod \limits_{j=m+1}^d \big| \lambda_k c_j^k \big|} \\
&&~~~~\xi(\dot{r},\ddot{r})~d(\dot{r},\ddot{r}) \\
&&~\\
&=& \int \limits_{\R^d} \int \limits_{\R^{d-m}} \hat f_k^{3,4,5,6,7,8}\Big(\dot{s}-\dot{r},\ddot{s}-\ddot{r},-\frac{\lambda_k \dot{c}^k}{2}(\dot{s}+\dot{r}),-\frac{\lambda_k \ddot{c}^k}{2}(\ddot{s}+\ddot{r}),\rho_k,\lambda_k,0,0\Big) \\
&&~~~~~~~~~~~\quer{\ddot{\eta}}_{k,0}(\tilde{y})\ddot{\eta}_{k,0}(\tilde{y})~d\tilde{y}~\xi(\dot{r},\ddot{r})~d(\dot{r},\ddot{r}) \\ 
&&~\\
&-& \int \limits_{\R^{d}} \int \limits_{\R^{d-m}} \hat{f}_{\infty}^{3,4,5,6,7,8}\Big(\dot{s}-\dot{r},\ddot{s}-\ddot{r},-\frac{\la}{2} \dot{c}(\dot{s}+\dot{r}),\la_k \ddot{c}^k(\tilde{y}-\ddot{s}),\rho,\lambda,0,0 \Big) \\ 
&&~~~~~~~~~~~\quer{\ddot{\eta}}_{k,0}(\tilde{y}+\ddot{r}-\ddot{s}) \ddot{\eta}_{k,0}(\tilde{y})~d\tilde{y} ~\xi(\dot{r},\ddot{r})~d(\dot{r},\ddot{r}).
\end{eqnarray*}
~\\
Similar as for the second case, functions $q_k$, $u_k$, $v_k$, $o_k$, $n_k$ and $w_k$ are going to be defined in order to divide the above integrals into six parts:
\begin{eqnarray*}
q_k(s, \tilde{y})&:=&\int \limits_{\R^{d}} \xi(\dot{r},\ddot{r})\quer{\ddot{\eta}}_{k,0}(\tilde{y}+\ddot{r}-\ddot{s}) \\
&&~~~~\bigg(\hat{f}_k^{3,4,5,6,7,8}\Big(\dot{s}-\dot{r},\ddot{s}-\ddot{r},-\frac{\lambda_k \dot{c}^k}{2}(\dot{s}+\dot{r}),-\frac{\lambda_k \ddot{c}^k}{2}(\ddot{s}+\ddot{r}),\textcolor{red}{\rho_k},\lambda_k,0,0\Big) \\ 
&&~~~~~~~-\hat{f}^{3,4,5,6,7,8}_k \Big(\dot{s}-\dot{r},\ddot{s}-\ddot{r},-\frac{\lambda_k \dot{c}^k}{2}(\dot{s}+\dot{r}),-\frac{\lambda_k \ddot{c}^k}{2}(\ddot{s}+\ddot{r}),\textcolor{red}{\rho},\lambda_k,0,0 \Big)\bigg)d(\dot{r},\ddot{r}),
\end{eqnarray*}
\begin{eqnarray*}
u_k(s, \tilde{y})&:=&\int \limits_{\R^{d}} \xi(\dot{r},\ddot{r})\quer{\ddot{\eta}}_{k,0}(\tilde{y}+\ddot{r}-\ddot{s})\\
&&~~~~\bigg(\hat{f}_k^{3,4,5,6,7,8}\Big(\dot{s}-\dot{r},\ddot{s}-\ddot{r},-\frac{\lambda_k \dot{c}^k}{2}(\dot{s}+\dot{r}),-\frac{\lambda_k \ddot{c}^k}{2}(\ddot{s}+\ddot{r}),\rho,\textcolor{red}{\lambda_k},0,0 \Big) \\ 
&&~~~~~~~-\hat{f}_k^{3,4,5,6,7,8}\Big(\dot{s}-\dot{r},\ddot{s}-\ddot{r},-\frac{\lambda_k \dot{c}^k}{2}(\dot{s}+\dot{r}),-\frac{\lambda_k \ddot{c}^k}{2}(\ddot{s}+\ddot{r}),\rho,\textcolor{red}{\lambda},0,0 \Big)\bigg)d(\dot{r},\ddot{r}),
\end{eqnarray*}
\begin{eqnarray*}
v_k(s, \tilde{y})&:=&\int \limits_{\R^{d}}\xi(\dot{r},\ddot{r})\quer{\ddot{\eta}}_{k,0}(\tilde{y}+\ddot{r}-\ddot{s})\\
&&~~~~\bigg(\hat{f}_k^{3,4,5,6,7,8}\Big(\dot{s}-\dot{r},\ddot{s}-\ddot{r},-\frac{\lambda_k \dot{c}^k}{2}(\dot{s}+\dot{r}),\textcolor{red}{-\frac{\lambda_k \ddot{c}^k}{2}(\ddot{s}+\ddot{r})},\rho,\lambda,0,0 \Big) \\ 
&&~~~~~~~-\hat{f}_k^{3,4,5,6,7,8}\Big(\dot{s}-\dot{r},\ddot{s}-\ddot{r},-\frac{\lambda_k \dot{c}^k}{2}(\dot{s}+\dot{r}),\textcolor{red}{\lambda_k \ddot{c}^k(\tilde{y}-\ddot{s})}, \rho,\lambda,0,0 \Big)\bigg)d(\dot{r},\ddot{r}),
\end{eqnarray*}
\begin{eqnarray*}
o_k(s, \tilde{y})&:=&\int \limits_{\R^{d}}\xi(\dot{r},\ddot{r})\quer{\ddot{\eta}}_{k,0}(\tilde{y}+\ddot{r}-\ddot{s})\\
&&~~~~\bigg(\hat{f}_k^{3,4,5,6,7,8}\Big(\dot{s}-\dot{r},\ddot{s}-\ddot{r},\textcolor{red}{-\frac{\lambda_k \dot{c}^k}{2}(\dot{s}+\dot{r})},\lambda_k \ddot{c}^k(\tilde{y}-\ddot{s}), \rho,\lambda,0,0 \Big) \\ 
&&~~~~~~~-\hat{f}_k^{3,4,5,6,7,8}\Big(\dot{s}-\dot{r},\ddot{s}-\ddot{r},\textcolor{red}{-\frac{\la}{2} \dot{c}(\dot{s}+\dot{r})},\la_k \ddot{c}^k(\tilde{y}-\ddot{s}),\rho,\lambda,0,0\Big)\bigg)d(\dot{r},\ddot{r}),
\end{eqnarray*}
\begin{eqnarray*}
n_k(s, \tilde{y})&:=&\int \limits_{\R^{d}}\xi(\dot{r},\ddot{r})\quer{\ddot{\eta}}_{k,0}(\tilde{y}+\ddot{r}-\ddot{s})\\
&&~~~~\bigg(\textcolor{red}{\hat{f}_{k}^{3,4,5,6,7,8}}\Big(\dot{s}-\dot{r},\ddot{s}-\ddot{r},-\frac{\la}{2} \dot{c}(\dot{s}+\dot{r}),\la_k \ddot{c}^k(\tilde{y}-\ddot{s}),\rho,\lambda,0,0\Big) \\ 
&&~~~~~~~-\textcolor{red}{\hat{f}_{\infty}^{3,4,5,6,7,8}}\Big(\dot{s}-\dot{r},\ddot{s}-\ddot{r},-\frac{\la}{2} \dot{c}(\dot{s}+\dot{r}),\la_k \ddot{c}^k(\tilde{y}-\ddot{s}),\rho,\lambda,0,0\Big)\bigg)d(\dot{r},\ddot{r})
\end{eqnarray*}
and
\begin{eqnarray*}
w_k(s)&:=&\int \limits_{\R^{d-m}} \int \limits_{\R^d} \xi(\dot{r},\ddot{r})\ddot{\eta}_{k,0}(\tilde{y})\Big(~\quer{\ddot{\eta}}_{k,0}(\tilde{y})-\quer{\ddot{\eta}}_{k,0}(\tilde{y}+\ddot{r}-\ddot{s})\Big) \\ 
&&~~~~~~~~~~\hat{f}_k^{3,4,5,6,7,8}\Big(\dot{s}-\dot{r},\ddot{s}-\ddot{r},-\frac{\lambda_k \dot{c}^k}{2}(\dot{s}+\dot{r}),-\frac{\lambda_k \ddot{c}^k}{2}(\ddot{s}+\ddot{r}),\rho_k,\lambda_k,0,0\Big)~d (\dot{r},\ddot{r})d\tilde{y}.
\end{eqnarray*}
Then,
\begin{eqnarray*}
\big(\pi_k(f)- \nu_k\big(p_{G/U}(f)\big)\big)\xi (s)
&=&\int \limits_{\R^{d-m}} q_k(s, \tilde{y}) \ddot{\eta}_{k,0}(\tilde{y})~d \tilde{y}~+~\int \limits_{\R^{d-m}} u_k(s, \tilde{y}) \ddot{\eta}_{k,0}(\tilde{y})~d \tilde{y} \\
&+&\int \limits_{\R^{d-m}} v_k(s, \tilde{y}) \ddot{\eta}_{k,0}(\tilde{y})~d \tilde{y}
~+~\int \limits_{\R^{d-m}} o_k(s, \tilde{y}) \ddot{\eta}_{k,0}(\tilde{y})~d \tilde{y} \\
&+&\int \limits_{\R^{d-m}} n_k(s, \tilde{y}) \ddot{\eta}_{k,0}(\tilde{y})~d \tilde{y}~+~w_k(s).
\end{eqnarray*}
~\\
As in the second case, to show that
$$\big\|\pi_k(f)- \nu_k\big(p_{G/U}(f)\big)\big\|_{op} \overset{k \to \infty}{\longrightarrow} 0,$$
one has to prove that there are $\kappa_k$, $\gamma_k$, $\delta_k$, $\tau_k$, $\omega_k$ and $\epsilon_k$ which are going to $0$ for $k \to \infty$, such that 
$$\|q_k\|_2\klgl\kappa_k \| \xi \|_2, ~~~\|u_k\|_2\klgl \gamma_k \| \xi \|_2, ~~~\|v_k\|_2 \klgl \delta_k \| \xi \|_2,~~~\|o_k\|_2 \klgl \tau_k \| \xi \|_2,$$
$$\|n_k\|_2 \klgl \omega_k \| \xi \|_2~~~\text{and}~~~\|w_k\|_2 \klgl \epsilon_k \| \xi \|_2.$$

The estimation of the functions $q_k$, $u_k$, $v_k$, $n_k$ and $w_k$ is very similar to their estimation in the second case and will thus be skipped. So, it just remains the estimation of $o_k$: \\
~\\
For this, first regard the last factor of the function $o_k$: 
\begin{eqnarray*}
&&\hat{f}_k^{3,4,5,6,7,8}\Big(\dot{s}-\dot{r},\ddot{s}-\ddot{r},-\frac{\lambda_k \dot{c}^k}{2}(\dot{s}+\dot{r}),\lambda_k \ddot{c}^k(\tilde{y}-\ddot{s}), \rho,\lambda,0,0 \Big) \\  
&-&\hat{f}_k^{3,4,5,6,7,8}\Big(\dot{s}-\dot{r},\ddot{s}-\ddot{r},-\frac{\la}{2} \dot{c}(\dot{s}+\dot{r}),\la_k \ddot{c}^k(\tilde{y}-\ddot{s}),\rho,\lambda,0,0 \Big) \\
&=& \Big(~\frac{1}{2}(\la \dot{c}-\lambda_k \dot{c}^k)(\dot{s}+\dot{r})\Big) \\
&\cdot&\int \limits_0^1 \partial_3 \hat{f}_k^{3,4,5,6,7,8}\Big(\dot{s}-\dot{r},\ddot{s}-\ddot{r},-\frac{\la}{2} \dot{c}(\dot{s}+\dot{r})+t\Big(~\frac{1}{2}(\la \dot{c}-\lambda_k \dot{c}^k)(\dot{s}+\dot{r})\Big),\lambda_k \ddot{c}^k(\tilde{y}-\ddot{s}), \rho,\lambda,0,0\Big)dt.
\end{eqnarray*}

Thus, there exists an on $f$ depending constant $C_1>0$ such that 
\begin{eqnarray*}
&&\Big|~\hat{f}_k^{3,4,5,6,7,8}\Big(\dot{s}-\dot{r},\ddot{s}-\ddot{r},-\frac{\lambda_k \dot{c}^k}{2}(\dot{s}+\dot{r}),\lambda_k \ddot{c}^k(\tilde{y}-\ddot{s}), \rho,\lambda,0,0\Big) \\
&&~~-\hat{f}_k^{3,4,5,6,7,8}\Big(\dot{s}-\dot{r},\ddot{s}-\ddot{r},-\frac{\la}{2} \dot{c}(\dot{s}+\dot{r}),\la_k \ddot{c}^k(\tilde{y}-\ddot{s}),\rho,\lambda,0,0\Big) \Big| \\
&\leq& \big\|\la \dot{c}-\lambda_k \dot{c}^k\big\| \|\dot{s}+\dot{r}\|~ \frac{C_1}{(1+\|\dot{s}+\dot{r}\|)^{2d+1}~(1+\|\ddot{s}-\ddot{r}\|)^{2d}}.
\end{eqnarray*}

Hence, one gets
\begin{eqnarray*}
&&\| o_k \|_2^2 \\
&&~\\
&=& \int \limits_{\R^{d+(d-m)}} |o_k(s,\tilde{y})|^2 d(s, \tilde{y}) \\
&\leq& \int \limits_{\R^{d+(d-m)}}\bigg(~\int \limits_{\R^d}|\xi(\dot{r},\ddot{r})| \big|\quer{\ddot{\eta}}_{k,0}(\tilde{y}+\ddot{r}-\ddot{s})\big| \big\|\la \dot{c}-\lambda_k \dot{c}^k \big\| \|\dot{s}+\dot{r}\|\\
&&~~~~~~~~~~~~\frac{C_1}{(1+\|\dot{s}+\dot{r}\|)^{2d+1}~(1+\|\ddot{s}-\ddot{r}\|)^{2d}}~ d(\dot{r},\ddot{r}) \bigg)^2 d(\dot{s},\ddot{s}, \tilde{y}) \\
&&~\\
&\underset{Schwarz}{\overset{Cauchy-}{\leq}}& C_1^2 \big\| \la \dot{c}-\lambda_k \dot{c}^k \big\|^2 \int \limits_{\R^{d+(d-m)}} \bigg( ~\int \limits_{\R^d} \frac{1}{(1+\|\dot{s}+\dot{r}\|)^{2d}~(1+\|\ddot{s}-\ddot{r}\|)^{2d}}~ d(\dot{r},\ddot{r})\bigg) \\
&& ~~\bigg(~ \int \limits_{\R^d} |\xi(\dot{r},\ddot{r})|^2 \big|\quer{\ddot{\eta}}_{k,0}(\tilde{y}+\ddot{r}-\ddot{s}) \big|^2 \frac{\|\dot{s}+\dot{r}\|^2}{(1+\|\dot{s}+\dot{r}\|)^{2d+2}~(1+\|\ddot{s}-\ddot{r}\|)^{2d}}~ d(\dot{r},\ddot{r}) \bigg) d(\dot{s},\ddot{s}, \tilde{y}) \\
&\leq&C'_1 \big\| \la \dot{c}-\lambda_k \dot{c}^k \big\|^2 \int \limits_{\R^{d+(d-m)}} \int \limits_{\R^d} |\xi(\dot{r},\ddot{r})|^2 \big|\quer{\ddot{\eta}}_{k,0}(\tilde{y}+\ddot{r}-\ddot{s})\big|^2 \\
&&~~~~~~~~~~~~~~~~~~~~~~~~~~~~~~\frac{1}{(1+\|\dot{s}+\dot{r}\|)^{2d}~(1+\|\ddot{s}-\ddot{r}\|)^{2d}}~ d(\dot{r},\ddot{r})d(\dot{s},\ddot{s}, \tilde{y}) \\
&\overset{\|\eta_{k,0}\|_2=1}{=}& C''_1 \big\| \la \dot{c}-\lambda_k \dot{c}^k \big\|^2 \| \xi \|_2^2,
\end{eqnarray*}
with matching constants $C'_1>0$ and $C''_1>0$, depending on $f$. Hence, $\tau_k:=\sqrt{C''_1 }~ \big\| \la \dot{c}-\lambda_k \dot{c}^k \big\|$ fulfills $\tau_k \overset{k \to \infty}{\longrightarrow} 0$ and
$$\| o_k\|_2~\leq~\tau_k \|\xi\|_2.$$
~\\
Thus, for those $f \in \S(\R^{2d+2+p}) \cong {\cal S}(G)$ whose Fourier transform in $[\g,\g]$ has a compact support,
$$\big\|\pi_k(f)- \nu_k\big(p_{G/U}(f))\big\|_{op}~=~\sup_{\substack{\xi \in L^2(\R^d)\\ \| \xi\|_2=1}}\big\|\big(\pi_k(f)- \nu_k\big(p_{G/U}(f)\big)\big)(\xi)\big\|_2
\overset{k \to \infty}{\longrightarrow}0.$$
As in the second case, because of the density in $C^*(G)$ of the set of Schwartz functions whose partial Fourier transform has a compact support, the claim follows for all $a \in C^*(G)$. \\
\qed
~\\
~\\
Now, the assertions for the sequence $\big(\pi_k^V\big)_k$ can be deduced: 
\subsubsection{Transition to $\big(\pi_k^V\big)_k$}
Again, because of the equivalence of the representations $\pi_k$ and $\pi_k^V$ for every $k \in \N$, there exist unitary intertwining operators
$$F_k: \H_{\pi_k^V} \cong L^2(\R^d) \to \H_{\pi_k} \cong L^2(\R^d)~~\text{with}~~ F_k \circ \pi_k^V(a)=\pi_k(a) \circ F_k~~~\forall~ a \in C^*(G).$$
With the injective $*$-homomorphism 
$$\tau:C_0 \big(\R^{2(d-m)},\K \big) \to C^*(G/U, \chi_{\ell}),~\F(a)_{|L((O_k)_k)} \mapsto p_{G/U}(a)$$
define $$\tilde{\nu}_k(\va):=F_k^* \circ (\nu_k \circ \tau)\big(\va_{| L((O_k)_k)}\big) \circ F_k~~~\forall~\va \in CB(S_{i-1}).$$ 
Then, like in the second case, $\tilde{\nu}_k$ complies with the demanded requirements and thus, the original representations $\big(\pi_k^V\big)_k$ fulfill Property 3(b).\\

Finally, one obtains the following result:
\begin{theorem}[Main result]
~\\
The $C^*$-algebra $C^*(G)$ of a connected real two-step nilpotent Lie group is isomorphic (under the Fourier transform) to the set of all operator fields $ \va $ defined over $ \wh G $ such that
\begin{enumerate}\label{}
\item $ \va(\ga)\in \K(\H_i) $ for every $i \in \{1,...,r\}$ and every $ \ga\in\GA_i$.
\item $\va \in l^{\infty}(\wh{G})$.
\item  The mappings $ \ga \mapsto \va(\ga) $ are norm continuous on the different sets $ \GA_i $.
\item  For any sequence $ (\ga_k)_{k\in\N} \subset \wh G$ going to infinity $ \lim \limits_{k\to\iy}\noop{\va(\ga_k)}=0 $.
\item  For $i \in \{1,...,r\}$ and any properly converging sequence  $\ol \ga=(\ga_k)_k\subset \GA_i   $ whose limit set $ L(\ol\ga)$ is contained in $S_{i-1} $  (taking a subsequence if necessary) and for the mappings \mbox{$\tilde{\nu}_k=\tilde{\nu}_{\ol\ga,k}: CB(S_{i-1})\to \B(\H_i)$} constructed in the preceding sections, one has
\begin{eqnarray*}
&\lim \limits_{k\to\iy }
\big\|
\va(\ga_k)-\tilde{\nu}_{k}\big(\va\res{S_{i-1}}\big)\big\|_{op}=0. 
\end{eqnarray*}
 \end{enumerate}
\end{theorem}

\section{Example: The free two-step nilpotent Lie groups of $3$ and $4$ \mbox{generators}}
In the case of the free two-step nilpotent Lie groups of $n=3$ and $n=4$ generators, the stabilizer of a linear functional $\ell$, the in Section 
\ref{cond 1,2,3a} constructed polarization $\p_{\ell}^V$, the coadjoint orbits, as well as the sets $S_i$ and $\Gamma_i$ can easily be calculated. \\
For $n=3$, there are coadjoint orbits of the dimensions $0$ and $2$ and for $n=4$, the dimensions $0$, $2$ and $4$ appear. \\
For the free two-step nilpotent Lie groups of $3$ generators, the third case regarded in the proof above does not appear: For this, one has to find a sequence of orbits $(O_k)_k$ whose limit set $L((O_k)_k)$ consists of orbits of the dimension strictly greater than $0$ but strictly smaller than $\dim{O_k}$. But as for $n=3$ only orbits of the dimensions $0$ and $2$ appear, such a sequence $(O_k)_k$ does not exist. However, for the free two-step nilpotent Lie groups of $4$ generators, this discussed third case exists. \\
For both $n=3$ and $n=4$, one can also see that the situation occurs where the polari\-zations $\p_{\ell}^V$ are discontinuous in $\ell$ on the set $\{\ell_{O'}|~O' \in (\g^*/G)_{2d}\}$. This shows the necessity of regarding the sets $\{\ell_{O'}|~O' \in (\g^*/G)_{(J,K)}\}$ instead. \\
Some calculations for the example of the free two-step nilpotent Lie groups of $3$ and $4$ generators can be found in the doctoral thesis of R.Lahiani (see \cite{lahi}).

\section{Appendix}
\begin{lemma}\label{lemma}
~\\
Let $V$ be a finite-dimensional euclidean vector space and $S$ an invertible, skew-symmetric endomorphism. Then $V$ can be decomposed into an orthogonal direct sum of two-dimensional $S$-invariant subspaces. 
\end{lemma}

Proof:\\
$S$ extends to a complex  endomorphism  $ S_\C $ on the
complexification $V_\C$ of $V$, which has purely
imaginary eigenvalues. \\
If $ i\la\in i\R $ is an eigenvalue, then also $
-i\la $ is a spectral element. Denote by $ E_{i\la} $ the corresponding
eigenspace. These eigenspaces are orthogonal to each other with respect to the
Hilbert space structure of $
V_\C $ coming from the euclidean scalar product $ \langle \cdot,\cdot\rangle $ on $V$.

Let for $ i\la  $ in the spectrum of $ S_\C $
\begin{eqnarray*}
V^\la:=(E_{i\la}+E_{-i\la})\cap V.
\end{eqnarray*}
If $ \la\ne 0 $, $ dim \big(V^\la \big) $ is even and $V^\la $ is $ S $-invariant and orthogonal to $V^{\la'} $, whenever $ |\la |\ne|\la'| $: \\
Indeed, one then has for $ x\in V^\la, x'\in V^{\la'} $ that
\begin{eqnarray*}
&x+iy\in E_{i\la}~~~\text{and}~~~ x-iy\in E_{-i\la} \textrm{ for some }y\in V~~~\text{as well as}\\
 & x'+iy'\in E_{i\la'}~~~\text{and}~~~ x'-iy'\in E_{-i\la'} \textrm{ for some }y'\in V.
\end{eqnarray*}
Therefore,
$$\langle{x+iy},{x'+iy'}\rangle=0~~~\text{and}~~~\langle{x-iy},{x'+iy'}\rangle=0.$$
Thus, one has
$$\langle{x},{x'+iy'}\rangle=0~~~\text{and hence}~~~\langle{x},{x'}\rangle=0.$$

Suppose that $dim \big(V^\la \big)>2$, choose a vector $ x\in V^\la $ of length 1 and let $y=S(x)$. Since $S_\C^2=-\la^2 \textrm{Id}$, both on $E_{i\la}$ and on $E_{-i\la}$,
 $$S(y)~=~S^2(x)~=~-\la^2 x.$$
This shows that $ W_1^\la:= \text{span} \{x,y\} $ is an $S$-invariant subspace of $
V^\la $. If $ V_1^\la $ denotes the orthogonal complement of $ W_1^{\la} $ in $ V^\la $,
then $ V_1^\la $ is $ S $-invariant, since $S^t=-S$. \\
In this way one can find a decomposition of $ V^\la $ into an orthogonal direct sum of
two-dimensional $ S $-invariant subspaces  $ W_j^\la $ and by summing up over the eigenvalues, one obtains the required decomposition of $V$. 

\section{Acknowledgements}
This work is supported by the Fonds National de la Recherche, Luxem\-bourg (Project \mbox{Code 3964572}). 

\end{document}